\documentclass[a4paper,11pt]{article}
\usepackage{amsmath,amsthm}
\usepackage{amsfonts}
\usepackage{amssymb,dsfont,hyperref,graphicx}
\usepackage{caption, subcaption, color,commath}
\usepackage{cases, xcolor}
\usepackage[left = 3cm, right = 3cm, top = 3cm , bottom = 3cm]{geometry}

\allowdisplaybreaks[1]
\newtheorem{theorem}{Theorem}
\newtheorem{proposition}{Proposition}

\newtheorem{lemma}{Lemma}
\newtheorem{definition}{Definition}
\newtheorem{remark}{Remark}

\newcommand{\kv}{\mathbb{E}}
\newcommand{\id}{\mathds{1}}
\newcommand{\xs}{\mathbb{P}}
\newcommand{\ps}{\mathbb{V}ar}

\newcommand{\rhoa}{\mathbb{R}}
\newcommand{\rhoad}{\mathbb{R}_+}
\newcommand{\nhoa}{\mathbb{N}^*}
\newcommand{\ehoa}{\mathcal{E}}
\newcommand{\hhoa}{\mathcal{H}}
\newcommand{\ahoa}{\mathcal{A}}
\newcommand{\mhoa}{\mathcal{M}}

\newcommand{\usup}[1]{\,\underset{#1}{\sup}\,}
\newcommand{\Ftr}{\mathcal{F}}
\newcommand{\pt}{\partial}

\newcommand{\chuan}[2]{\|#1\|_{#2}}
\newcommand{\ch}[1]{\|#1\|_2}
\newcommand{\vc}{\infty}

\newcommand{\kepa}[1]{\langle #1,1 \rangle}

\DeclareMathOperator{\argmin}{argmin\,}

\DeclareMathOperator{\Geom}{Geom}

\newcommand{\lb}{\ell}
\newcommand{\htl}{\hat{h}_{\ell}}

\newcommand{\hl}{\hat{h}_\ell}
\newcommand{\hlp}{\hat{h}_{\ell'}}
\newcommand{\hlml}{\hat{h}_{\hat\ell,\ell}}
\newcommand{\hllp}{\hat{h}_{\ell,\ell'}}
\newcommand{\Kl}{K_\ell(\gamma - \Gamma^1_i)}
\newcommand{\Kll}{K_\ell(\gamma - \Gamma^1_1)}
\newcommand{\red}[1]{\textcolor{red}{#1}}

\DeclareMathOperator{\Beta}{Beta}

\title{Estimating the Division Kernel of a Size-Structured Population}
\author{Van Ha Hoang \footnote{Laboratoire Paul Painlev\'e UMR CNRS 8524, Universit\'e Lille 1.\newline
\hspace*{0.2in}\textit{Email address}: van-ha.hoang@math.univ-lille1.fr}}
\date{\today}
\begin{document}
\maketitle

\begin{abstract}
We consider a size-structured population describing the cell divisions. The cell population is described by an empirical measure and we observe the divisions 
in the continuous time interval $[0,T]$. We address here the problem of estimating the division kernel $h$ (or fragmentation kernel) in case of complete data. 
An adaptive estimator of $h$ is constructed based on a kernel function $K$ with a fully data-driven bandwidth selection method. We obtain an oracle inequality and an exponential convergence rate, for which optimality is considered.
\end{abstract}

\section{Introduction}
Models for populations of dividing cells possibly differentiated by covariates such as size have made the subject of an abundant literature in recent years (starting from Athreya and Ney \cite{Athreya70}, 
Harris \cite{Harris63}, Jagers \cite{Jagers69}...) Covariates termed as `size' are variables that grow deterministically with time (such as volume, length, level of certain proteins, DNA content, \textit{etc}.) Such models of structured populations provide descriptions for the evolution of the size distribution, which can be interesting for applications. For instance, in the spirit of Stewart \textit{et al.} \cite{Stewart05}, we can imagine that each cell contains some toxicities whose quantity plays the role of the size. The asymmetric divisions of the cells, where one daughter contains more toxicity than the other, can lead under some conditions to the purge of the toxicity in the population by concentrating it into few lineages. These results are linked with the concept of aging for cell lineage. This concept has been tackled in many papers (e.g. Ackermann \textit{et al.} \cite{Akermann03}, Aguilaniu \textit{et al.} \cite{Aguilaniu03},  C-Y. Lai \textit{et al.} \cite{Lai02}, Evans and Steinsaltz \cite{Evans07}, Moseley \cite{Moseley13}...).  

\par Here we consider a stochastic individual-based model of size-structured population 
in continuous time, where individuals are cells undergoing asymmetric binary divisions and whose size is the quantity of toxicity  they contain. A cell containing a toxicity $x\in\rhoa_+$ divides at a rate $R>0$. The toxicity grows inside the cell with rate $\alpha>0$. When a cell divides, a random fraction $\Gamma \in [0,1]$ of the toxicity goes in the first daughter cell and $1-\Gamma$ in the second one. If $\Gamma = \frac 12$, the daughters are the same with toxicity $\frac x2$. We assume that $\Gamma$ has a symmetric distribution on $[0,1]$ with a density $h$ with respect to Lebesgue measure 
such that $\xs(\Gamma=0) = \xs(\Gamma = 1) = 0$. If $h$ is piked at $1/2$ (\textit{i.e.} $\Gamma \simeq 1/2$), then both daughters contain the same toxicity, \textit{i.e.} the half of their mother's toxicity. The more $h$ puts weight in the neighbourhood of $0$ and $1$, the more asymmetric the divisions are, with one daughter having little toxicity and the other a toxicity close to its mother's one. If we consider that having a lot of toxicity is a kind of senescence, then, the kurtosis of $h$ provides indication on aging phenomena  (see \cite{Lindner15}).

\par Modifications of this model to account for more complex phenomena have been considered in other papers. Bansaye and Tran \cite{Bansaye11}, Cloez \cite{Cloez11} or Tran \cite{Tran08} consider non-constant division and growth rates. Robert \textit{et al.} \cite{LydiaRobert14} studies whether divisions can occur only when a size threshold is reached. Our purpose here is to estimate the density $h$ ruling the divisions, and we stick to constant rates $R$ and $\alpha$ for the sake of simplicity. 
Notice that several similar models for binary cell division in discrete time also exist in the literature and have motivated statistical question as here, see for instance Bansaye et al. \cite{Bansaye08,Bansaye13}, Bercu \textit{et al.} \cite{Bercu09}, Bitseki Penda \cite{BitsekiPenda15}, Delmas and Marsalle \cite{Delmas10} or Guyon \cite{Guyon07}.

\begin{figure}
\centering
 \includegraphics[scale = 0.6]{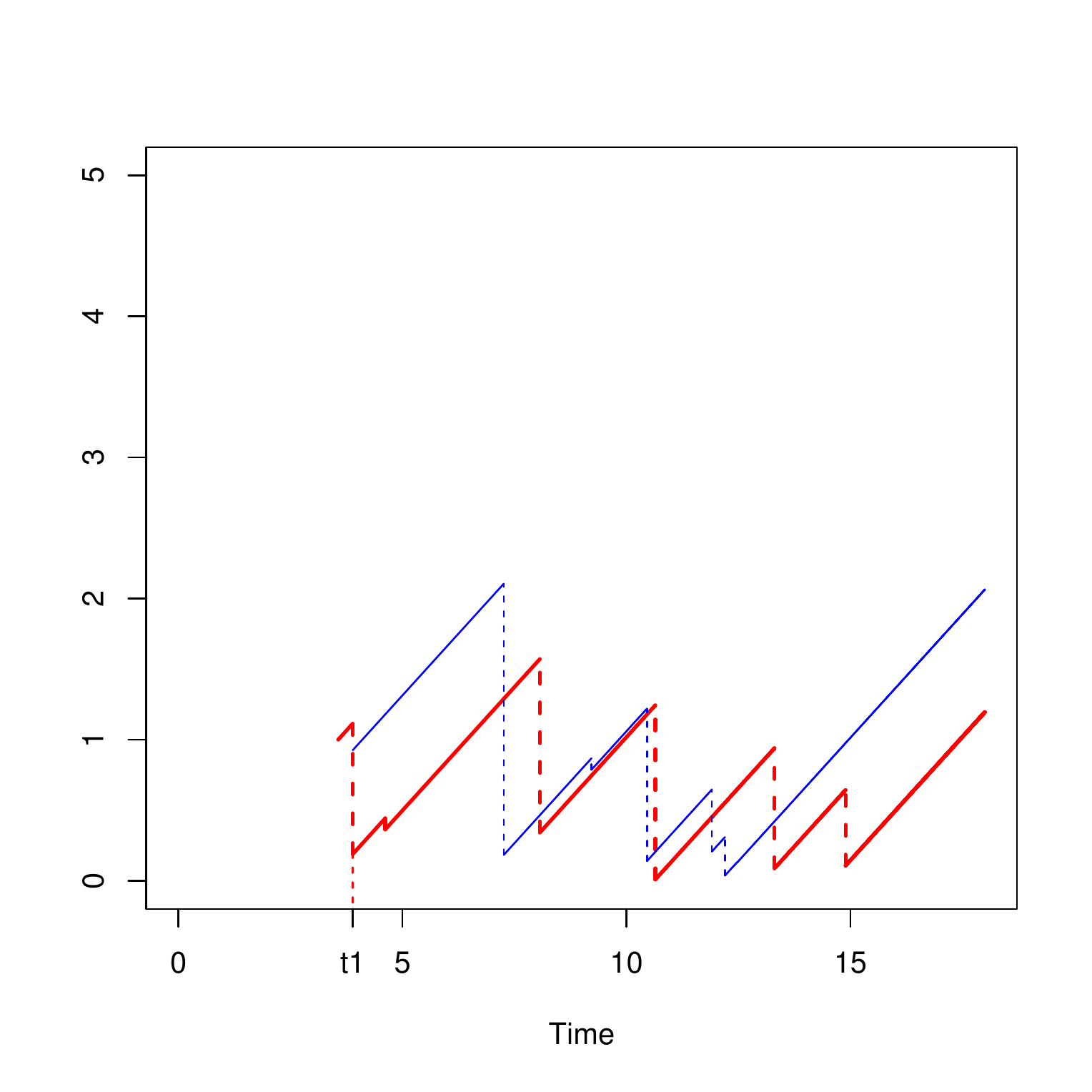}
\caption{{\small \textit{Trajectories of two daughter cells after a division, separating after the first division at time $t_1$.}} \label{fig2}}
\end{figure}

\par Individual-based models provide a natural framework for statistical estimation. Estimation of the division rate is, for instance, the subject of Doumic 
et al. \cite{DoumicHoffmann12, Doumic12} and Hoffmann and Olivier \cite{Hoffmann14}. Here, the density $h$ is the kernel division that we want to estimate. Assuming that we observe the divisions of cells in continuous time on the interval $[0,T]$, with $T>0$, we propose an adaptive kernel estimator $\hat{h}$ of $h$ 
for which we obtain an oracle inequality in Theorem \ref{th1}. The construction of $\hat{h}$ is detailed in the sequel. From oracle inequality we can infer adaptive exponential rates of convergence with respect to $T$ depending on $\beta$ the smoothness of the density.  Most of the time, nonparametric rates are of the form $n^{-\frac{2\beta}{2\beta+1}}$  (see for instance Tsybakov \cite{Tsybakov04}) and exponential rates are not often encountered in the literature. The exponential rates are due to binary splitting, the number of cells \textit{i.e} the sample size increases exponentially in $\exp(RT)$ (see Section 2.3). By comparison, in \cite{Hoffmann14} Hoffmann and Olivier obtain a similar rate of convergence $\exp\left(-\lambda_B\frac{\varsigma}{2\varsigma+1}T \right)$ of the kernel estimator of their division rate $B(x)$, where $\lambda_B$ is the Malthus parameter and $\varsigma>0$ is the smoothness of $B(x)$. However, their estimator $\hat{B}_T$ of $B$ is not adaptive since the choice of their optimal bandwidth still depends on $\varsigma$. Our estimator is adaptive with an ``optimal" bandwidth chosen from a data-driven method. We derive upper bounds and lower bounds for asymptotic minimax risks on H\"older classes and show that they coincide. Hence, the rate of convergence of our estimator $\hat h $ proves to be optimal in the minimax sense on the  H\"older classes. 

\par This paper is organized as follows. In Section 2, we introduce a stochastic differential equation driven by a Poisson point measure to describe the population of cells. Then, we construct the estimator of $h$ and obtain upper and lower bounds for the MISE (Mean Integrated Squared Error). Our main results are stated in Theorems \ref{th:upper-bound} and \ref{th:lower-bound}. Numerical results and discussions about aging effect are presented in Section 3. The main proofs are shown in Section 4.
\bigskip 

\textbf{Notation} \quad We introduce some notations used in the sequel. 

Hereafter, $\|\cdot\|_1$ and $\|\cdot\|_2$ denote the $\mathbb{L}^1$ and $\mathbb{L}^2$ norms on $\rhoa$ with respect to Lebesgue measure:
\[
\chuan{f}{1} = \int_\rhoa |f(\gamma)|d\gamma, \quad \ch{f} = \left(\int_\rhoa |f(\gamma)|^2d\gamma\right)^{1/2}.
\]
The $\mathbb{L}^\vc$ norm is defined by
\[
\chuan{f}{\vc} = \underset{\gamma\in (0,1)}{\sup}|f(\gamma)|.
\]
Finally, $f\star g$ denotes the convolution of two functions $f$ and $g$ defined by
\[
f\star g(\gamma) = \int_\rhoa f(u)g(\gamma -u)du.
\] 

\section{Microscopic model and kernel estimator of $h$}
\subsection{The model}
We recall the Ulam-Harris-Neveu notation used to describe the genealogical tree. The first cell is labelled by $\emptyset$ and when the cell $i$ divides, 
the two descendants are labelled by $i0$ and $i1$. The set of labels is
\begin{equation}\label{eq:Ulam-Harris-Neveu_notation}
J = \left\{\emptyset\right\} \cup \bigcup_{m=1}^{\infty}\left\{0,1\right\}^m.
\end{equation}
We denote $V_t$ the set of cells alive at time $t$, and $V_t\subset J$. 

\par Let $\mathcal{M}_F(\rhoa_+)$ be the space of finite measures on $\rhoad$ embedded with the topology of weak convergence and $X_t^i$ be the quantity of toxicity in the cell $i$ at time $t$, we describe the population of cells at time $t$ by a random point measure in $\mathcal{M}_F(\rhoa_+)$: 
\begin{equation}\label{eq:empircal-measure}
 Z_t(dx) = \sum_{i=1}^{N_t} \delta_{X^i_t}(dx),\quad\text{ where }\quad 
 N_t =\langle Z_t,1\rangle = \int_{\rhoad} Z_t(dx)
\end{equation}
is the number of individuals living at time $t$.
For a measure $\mu\in \mathcal{M}_F(\rhoa_+)$ and a positive function $f$, we use the notation $\langle \mu, f \rangle = \int_{\rhoad} fd\mu$.

\par Along branches of the genealogical tree, the toxicity $(X_t, t\ge 0)$ satisfies
\begin{equation}\label{eq:model-toxicity}
dX_t = \alpha dt,
\end{equation}
with $X_0 = x_0$. When the cells divide, the toxicity is shared between the daughter cells. This is described by the following stochastic differential equation (SDE). 

\par Let $Z_0 \in \mathcal{M}_F(\rhoa_+)$ be an initial condition such that
\begin{equation}\label{eq:initial-condition}
\kv(\kepa{Z_0})<+\vc,
\end{equation}
and let $Q(ds,di,d\gamma)$ be a Poisson point measure on $\rhoad\times \ehoa:=\rhoad\times J \times [0,1]$ with intensity $q(ds,di,d\gamma)= R\,ds\,n(di)H(d\gamma)$. $n(di)$ is the counting measure on $J$ and 
$ds$ is Lebesgue measure on $\rhoa_+$. We denote $\{\mathcal{F}_t \}_{t\ge 0}$ the canonical filtration associated with the Poisson point measure and the 
initial condition. The stochastic process $(Z_t)_{t\ge 0}$ can be described by a SDE as follows.

\begin{definition}\label{def:SDE-model}
For every test function $f_t(x) = f(x,t) \in \mathcal{C}^{1,1}_b(\rhoa_+\times \rhoa_+,\rhoa)$ (bounded of class $\mathcal{C}^1$ in $t$ and $x$ with bounded 
derivatives), the population of cells is described by:
\begin{align}
 &\langle Z_t,f_t\rangle = \langle Z_0, f_0\rangle + \int_0^t\int_{\rhoad} \big(\pt_s f_s(x) + \alpha\pt_x f_s(x)\big) Z_s(dx)ds \notag \\
&+ \int_0^t\int_{\ehoa} \mathds{1}_{\{i \le N_{s-}\}}\Big[f_s\left(\gamma X^i_{s-}\right) 
+ f_s\left((1-\gamma)X^i_{s-}\right) - f_s\left(X^i_{s-}\right)\Big]Q(ds,di,d\gamma). \label{eq:SDE-model}
\end{align}
\end{definition}
The second term in the right hand side of (\ref{eq:SDE-model}) corresponds to the growth of toxicities in the cells and the third term
gives a description of cell divisions where the sharing of toxicity into two daughter cells depends on the random fraction $\Gamma$. 

\par We now state some properties of $N_t$ that are useful in the sequel.

\begin{proposition}\label{prop:N_T}Let $T>0$, and assume the initial condition $N_0$, the number of mother cells at time $t=0$, is deterministic, for the sake of simplicity. We have
\begin{itemize}
\item[i)] Let $T_j$ be the $j^\text{th}$ jump time. Then:
\begin{equation}\label{eq:limit-N_T}
\lim_{j\to +\vc} T_j = +\vc \text{ and } \lim_{T\to+\vc} N_T = +\vc\quad\text{(a.s)}.
\end{equation}

\item[ii)] $N_T$ is distributed according to a negative binomial distribution, denoted  as $\mathcal{NB}(N_0,e^{-RT})$. Its probability mass function is then 
\begin{equation}\label{eq:dist-N_T-nb}
\xs\left(N_T = n \right) = \binom{n-1}{n-N_0} \big(e^{-RT}\big)^{N_0}\big(1 - e^{-RT}\big)^{n-N_0},
\end{equation}
for $n \ge N_0$. When $N_0 = 1$, $N_T$ has a geometric distribution
\begin{equation}\label{eq:dist-N_T-geom}
\xs\left(N_T = n \right) = e^{-RT}\left(1 - e^{-RT}\right)^{n-1}.
\end{equation}
Consequently, we have
\begin{equation}\label{eq:expectation-N_T}
\kv\big[N_T\big] = N_0e^{RT}.
\end{equation}
\item[iii)] When $N_0 = 1$:
\begin{equation}\label{eq:inverse-N_T-geom}
\kv\left[\frac{1}{N_T} \right] = \frac{RTe^{-RT}}{1-e^{-RT}}.
\end{equation}
When $N_0>1$, we have:
\begin{align}
\kv\left[\frac{1}{N_T} \right] &=
\left(\frac{e^{-RT}}{1-e^{-RT}} \right)^{N_0}(-1)^{N_0-1}\left(\sum_{k=1}^{N_0-1}\binom{N_0-1}{k} \frac{(-1)^ke^{kRT}}{k} + RT \right). \label{eq:inverse-N_T-nb} 
\end{align}
\item[iv)] Furthermore, when $N_0 >1$, we have
\begin{equation}\label{eq:bound-inverse-N_T-nb}
\frac{e^{-RT}}{N_0} \le \kv\left[\frac{1}{N_T} \right] \le \frac{e^{-RT}}{N_0-1}.
\end{equation}
\end{itemize}
\end{proposition}
The proof of Proposition \ref{prop:N_T} is presented in Section 4.

\subsection{Influence of age}\label{sec:mean-age}
In this section, we study the aging effect via the mean age which is defined as follows.

\begin{definition}\label{def:mean-age}
The mean age of the cell population up to time $t\in\rhoad$ is defined by:
\begin{equation}\label{eq:mean-age}
\bar{X}_t = \frac{1}{N_t}\sum_{i=1}^{N_t} X^i_t = \frac{\langle Z_t, f \rangle}{N_t}, 
\end{equation}
where $f(x) = x$.
\end{definition}

\par Following the work of Bansaye {et al.} \cite{Bansaye11a}, we note that the long time behavior of the mean age is related to the law of an auxiliary process $Y$ started at $Y_0 = \frac{X_0}{N_0}$ with infinitesimal generator characterized for all $f\in \mathcal{C}^{1,1}_b(\rhoa_+,\rhoa)$ by
\begin{equation}\label{eq:infinitesimal-generator}
 Af(x) = \alpha f'(x) + 2R\int_0^1 \left(f(\gamma x) - f(x) \right)h(\gamma) d\gamma.
\end{equation}

\par The empirical distribution $\frac{1}{N_t}\sum_{i=1}^{N_t}\delta_{X_t^i}$ gives the law of the path of a particle chosen at random at time $t$. Heuristically, the distribution of $Y$ restricted to $[0,t]$ approximates this distribution. Hence, this explains the coefficient $2$ which is a size-biased phenomenon, \textit{i.e.} when one chooses a cell in the population at time $t$, a cell belonging to a branch with more descendants is more likely to be chosen.

\begin{lemma}\label{lem1}
Let $Y$ be the auxiliary process with infinitesimal generator \eqref{eq:infinitesimal-generator}, for $t\in\rhoad$,
\begin{equation}\label{eq:SDE-mean-age}
Y_t = \left(Y_0 - \frac{\alpha}{R} \right)e^{-Rt} + \frac{\alpha}{R} + \int_0^t e^{-R(t-s)}dU_s.
\end{equation}
where $U_t$ is a square-integrable martingale.

\par Consequently, we have
\begin{equation}\label{eq:ave-mean-age}
\kv\left[Y_t \right] = \left(Y_0 - \frac{\alpha}{R}\right)e^{-Rt} + \frac{\alpha}{R},
\end{equation}
and
\begin{equation}\label{eq:limit-ave-mean-age}
\underset{t\to\vc}{\lim} \kv\left[Y_t \right] = \frac{\alpha}{R}.
\end{equation}
\end{lemma}

\par We will show that the auxiliary process $Y$ satisfies ergodic properties (see Section \ref{proof:th-mean-age}) which entails the following theorem.

\begin{theorem}\label{th:mean-age}
Assume that there exists $\underline{h}>0$ such that for all $\gamma \in (0,1)$, $h(\gamma) \ge \underline{h}$. Then
\begin{equation}\label{eq:limit-mean-age}
\underset{t\to +\vc}{\lim} \bar{X}_t = \underset{t\to +\vc}{\lim} \kv(Y_t) = \frac{\alpha}{R}.
\end{equation}
\end{theorem} 

\par Theorem \ref{th:mean-age} is a consequence of the ergodic properties of $Y$, of Theorem 4.2 in Bansaye {et al.} \cite{Bansaye11a} and of Lemma \ref{lem1}. It shows that the average of the mean age tends to the constant $\alpha/R$ when the time $t$ is large. Simulations in Section 3 illustrate the results. The proofs of Lemma \ref{lem1} and Theorem \ref{th:mean-age} are presented in Section \ref{proof:lem-mean-age} and Section \ref{proof:th-mean-age}.

\begin{remark}
When the population is large, we are interested in studying the asymptotic behavior of the random point measure. As in Doumic {et al.} \cite{Doumic12}, 
we can show that our stochastic model is approximated by a growth-fragmentation partial differential equation. This problem is a work in progress. 
\end{remark}

\subsection{Estimation of the division kernel}
\noindent\textbf{Data and construction of the estimator}\\[6pt]
Suppose that we observe the evolution of the cell population in a given time interval $[0,T]$. At the $i^{\text{th}}$ division time $t_i$, let us denote $j_i$ 
the individual who splits into two daughters $X_{t_i}^{j_i0}$ and $X_{t_i}^{j_i1}$ and define
 \[
  \Gamma^0_i = \frac{X_{t_i}^{j_i0}}{X_{t_i-}^{j_i}} \quad\text{ and }\quad \Gamma^1_i = \frac{X_{t_i}^{j_i1}}{X_{t_i-}^{j_i}},
 \]
the random fractions that go into the daughter cells, with the convention $\frac{0}{0} =0$.

\par $\Gamma^0_i$ and $\Gamma^1_i$ are exchangeable with $\Gamma^0_i + \Gamma^1_i =1$, $\Gamma^0_i$ and $\Gamma^1_i$ are thus not independent 
but the couples $(\Gamma^0_i,\Gamma^1_i)_{i\in\nhoa}$ are independent and identically distributed with distribution $(\Gamma^0,\Gamma^1)$ where 
$\Gamma^1\sim H(d\gamma)$ and $\Gamma^0 = 1 - \Gamma^1$. 

\par Since $h$ is a density function, it is natural to use a kernel method. We define an estimator $\hat{h}_\ell$ of $h$ based on the data $(\Gamma^0_i,\Gamma^1_i)_{i\in\nhoa}$ as follows.

\begin{definition}\label{def:estimator-h}
Let $K:\rhoa \longrightarrow \rhoa$ is an integrable function such that
\[
\int_\rhoa K(x)dx = 1 \;\text{ and }\; \int_\rhoa K^2(x)dx <\infty.
\] 
Let $M_T$ be the random number of divisions in the time interval $[0,T]$ and assume that $M_T > 0$. For all $\gamma\in (0,1)$, define
\begin{equation}\label{eq:estimator-h}
 \hat{h}_{\ell}(\gamma) = \frac{1}{M_T}\sum_{i=1}^{M_T} K_\ell(\gamma - \Gamma^1_i),
\end{equation}
where $K_\ell = \frac 1\ell K(\cdot/\ell)$, $\ell>0$ is the bandwidth to be chosen.
\end{definition}

\begin{remark}
Since $N_0\neq 0$, the number of random divisions $M_T$ is not equal to the number of individuals living at time $T$. Indeed, we have $M_T = N_T - N_0$. 
\end{remark}

In \eqref{eq:estimator-h}, $\hat{h}_\ell$ depends also on $T$. However, we omit $T$ for the sake of notation. The estimator $\hat{h}_{\ell}$ will satisfy the following properties.

\begin{proposition}\label{prop:mean-var-h}\quad
\begin{itemize}
\item[i)] The conditional expectation and conditional variance given $M_T$ of $\htl(\gamma)$ and variance $\htl(\gamma)$ are:
\begin{align}
\kv\big[\hat{h}_{\ell}(\gamma)|M_T\big] &= K_\ell\star h(\gamma) \,\text{ and }\, \kv\big[\hat{h}_{\ell}(\gamma) \big] = K_\ell\star h(\gamma) \label{eq:mean-h},\\
 \ps\big[\hat{h}_{\ell}(\gamma)\big|M_T] &= \frac{1}{M_T}\ps\left[K_\ell(\gamma - \Gamma^1_1)\right], \label{eq:var-cond-h} \\
  \ps\big[\hat{h}_{\ell}(\gamma)\big] &= \kv\big[\frac{1}{M_T}\big]\ps\left[K_\ell(\gamma - \Gamma^1_1)\right] \label{eq:var-h}.
\end{align}
Consequently, we have $\kv\big[\hat{h}_{\ell}(\gamma)|M_T\big]  = \kv\big[\hat{h}_{\ell}(\gamma) \big]$.
\item[ii)] For all $\gamma\in (0,1)$,
\begin{equation}\label{eq:assymp-h}
 \lim_{T\to +\vc}\hat{h}_{\ell}(\gamma) = K_\ell\star h(\gamma)\quad\text{(a.s)}.
\end{equation}
\end{itemize}
\end{proposition}

\noindent\textbf{Adaptive estimation of $h$ by Goldenshluger and Lepski's (GL) method}\\[6pt]
Let $\hat{h}_{\ell}$ be the kernel estimator of $h$ as in Definition \ref{def:estimator-h}. We measure the performance of $\hat{h}_{\ell}$ via its $\mathbb{L}^2$-loss \textit{i.e}  the average $\mathbb{L}^2$ distance between $\hat h_\ell$ and $h$. 
The objective is to find a bandwidth which minimizes this $\mathbb{L}^2$-loss. Since $M_T$ is random, we first study the  $\mathbb{L}^2$-loss conditionally to $M_T$. 
\begin{proposition}\label{prop:bias-variance-decomposition}
The $\mathbb{L}^2$-loss of $\hat h_\ell$ given $M_T$ satisfies :
\begin{equation}\label{eq:bias-variance-decomposition}
\kv\Big[\ch{\htl - h}\big| M_T \Big] \le \ch{h - K_\ell\star h} + \frac{\ch{K}}{\sqrt{M_T\ell}}.
\end{equation}
\end{proposition}

In the right hand side of the risk decomposition \eqref{eq:bias-variance-decomposition} the first term  is a bias term. Hence it decreases when $\ell\to 0$ whereas the second term which is a variance term increases when $\ell\to 0$. The best choice of $\ell$ should minimize this bias-variance trade-off. Thus, from a finite family of bandwidths $H$, the best  bandwidth $\bar\ell$ would be 
\begin{equation}\label{eq:oracle-band1}
\bar\ell := \underset{\ell\in H}{\argmin}\Big\{\ch{h - K_\ell\star h} + \frac{\ch{K}}{\sqrt{M_T\ell}} \Big\}.
\end{equation}

The bandwidth $\bar{\ell}$ is called "the oracle bandwidth" since it depends on $h$ which is unknown and then it cannot be used in practice. Since the oracle bandwidth minimizes a bias variance trade-off, we need to find an estimation  for the bias-variance decomposition of $\hat h_\ell$. Goldenshluger and Lepski \cite{GL11} developed a fully data-driven bandwidth selection method (GL method). The main idea of this method is based on an estimate of  the bias term by looking at several estimators. In a similar fashion, Doumic {et al.} \cite{Doumic12} and Reynaud-Bouret {et al.} \cite{R-B14} have used  this method. To apply the GL method, we set for any $\ell$, $\ell'\in H$:
\[
\hat h_{\lb,\lb'}:= \frac{1}{M_T}\sum_{i=1}^{M_T} \big(K_\lb\star K_{\lb'}\big)(\gamma - \Gamma^1_i)
= \left(K_\ell\star\hlp \right)(\gamma).
\]
Finally, the adaptive bandwidth and the estimator of $h$ are selected as follows:
\begin{definition}\label{def:GL-bw}
Given $\epsilon >0$ and setting $\chi:= (1+\epsilon)(1+\chuan{K}{1})$, we define
\begin{equation}\label{eq:data-driven-bandwidth}
\hat\lb := \underset{\lb\in H}{\argmin}\Big\{A(\lb) +  \frac{\chi\ch{K}}{\sqrt{M_T\lb}}\Big\},
\end{equation}
where, for any $\ell\in H$, 
\begin{equation}\label{eq:bias-estimator}
A(\lb):=\underset{\lb'\in H}{\sup}\Big\{\ch{\hat h_{\lb,\lb'} - \hat h_{\lb'}} 
- \frac{\chi\ch{K}}{\sqrt{M_T\lb'}}\Big\}_{+},
\end{equation}
Then, the estimator $\hat{h}$ is given by 
\begin{equation}\label{eq:adaptive-estimator-h}
\hat{h} := \hat{h}_{\hat\ell}.
\end{equation}
\end{definition}

An inspection of the proof of Theorem 2 shows that the term $ A(\lb)$ provides a control for the bias $\| h- K_\ell \star h\|_2$ up to the term $\| K\|_1$ (see (\ref{Aell}) and  (\ref{proof-th1-El}) in the proof of Theorem \ref{th1}, section 4). Since $A(\lb)$ depends only on $\hllp$ and $\hlp$, the estimator $\hat{h}$ can be computed in practice.   

\par We shall now state an oracle inequality which highlights the bias-variance decomposition of the MISE of $\hat h$. We recall that the MISE of $\hat h$ is the quantity $ \kv\Big[\ch{\hat h - h}^2 \Big] $.

\begin{theorem}\label{th1}
Let $T>0$ and assume that observations are taken on $[0,T]$. Let $N_0$ be the number of mother cells at the beginning of divisions and $M_T$ is the random number of divisions in $[0,T]$. Consider $H$ a countable subset of $\{\triangle^{-1}:\triangle=1,\ldots,\triangle_{\max}\}$ in which we choose the bandwidths and 
$\triangle_{\max}=\lfloor \delta M_T \rfloor$ for some $\delta>0$. Assume $h\in L^\vc([0,1])$ and let $\hat{h}$ be a kernel estimator defined with the kernel $K_{\hat\ell}$ where $\hat\ell$ is chosen by the GL method. Define
\begin{equation}\label{eq:rate}
\varrho(T)^{-1} = \begin{cases}
\dfrac{e^{-RT+\log(RT)}}{1-e^{-RT}},&\mbox{if } N_0 = 1, \\
e^{-RT},&\mbox{if } N_0 > 1.
\end{cases}
\end{equation}
For large $T$, the main term in $\varrho(T)$ is $e^{-RT}$ in any case. It is exactly the order of $\varrho(T)$ for $N_0>1$.
Then, given $\epsilon>0$
\begin{equation}\label{eq:oracle-ineq}
\kv\Big[\ch{\hat h - h}^2 \Big] \le 
C_1\underset{\ell \in H}{\inf}\left\{\ch{K_\ell\star h - h}^2 + \frac{\ch{K}^2}{\ell}\varrho(T)^{-1}\right\} + C_2\varrho(T)^{-1},
\end{equation}
where $C_1$ is a constant depending on $N_0$, $\chuan{K}{1}$ and $\epsilon$ and $C_2$ is a constant depending on $N_0$, $\delta$, $\epsilon$, $\chuan{K}{1}$, $\ch{K}$ and $\chuan{h}{\vc}$. 
\end{theorem}

The term $ \ch{K_\ell\star h - h}^2  $ is an approximation term, $  \frac{\ch{K}^2}{\ell}\varrho(T)^{-1}$ is a variance term and the last term  $ \varrho(T)^{-1}$ is asymptotically negligible. Hence the right hand side of the oracle inequality corresponds to a bias variance trade-off. 

\par We now establish upper and lower bounds for the MISE. The lower bound is obtained by perturbation methods (Theorem \ref{th:lower-bound}) and is valid for any estimator $\widehat{h}_T$ of $h$, thus indicating the optimal convergence rate. The upper bound is obtained in Theorem \ref{th:upper-bound} thanks to the key oracle inequality of Theorem \ref{th1}.

\par For the rate of convergence, it is necessary to assume that the density $h$ and the kernel function $K$ satisfy some regularity conditions introduced in the following definitions.

\begin{definition}\label{def:Holder-class}
Let $\beta>0$ and $L>0$. The H\"older class of smoothness $\beta$ and radius $L$ is defined by 
\begin{align*}
\hhoa(\beta,L) = \Big\{f:f\text{ has } &k = \lfloor\beta\rfloor \text{ derivatives and } \forall x,y\in\rhoa\\
&\big|f^{(k)}(y) - f^{(k)}(x)\big|\le L|x-y|^{\beta - k}\Big\}.
\end{align*}
\end{definition}

\begin{definition}\label{def:regularity-K}
Let $\beta^*>0$. An integrable function $K:\rhoa \to \rhoa$ is a kernel of order $\beta^*$ if 
\begin{itemize}
\item $\int K(x)dx = 1$,
\item $\int |x|^{\beta^*}|K(x)|dx <\vc$,
\item For $k = \lfloor \beta^* \rfloor$, $\forall 1\le j \le k$, $\int x^jK(x)dx=0$.
\end{itemize}
\end{definition}

Then, the following theorem gives  the rate of convergence of the adaptive estimator $\hat h$. 
\begin{theorem}\label{th:upper-bound}
Let $\beta^*>0$ and K be a kernel of order $\beta^*$. Let $\beta\in(0,\beta^*)$. Let $\hat\ell$ be the adaptive bandwidth defined in \eqref{eq:data-driven-bandwidth}. Then, for any $T>0$, the kernel estimator $\hat h$ satisfies
\begin{equation}\label{eq:upper-bound}
\underset{h\in \mathcal{H}(\beta,L)}{\sup}\kv\ch{\hat h - h}^2 \le C_3\varrho(T)^{-\frac{2\beta}{2\beta+1}},
\end{equation}
where $\varrho(T)^{-1}$ is defined in \eqref{eq:rate} and $C_3$ is a constant depending on $N_0$, $\delta$, $\epsilon$, $\chuan{K}{1}$, $\ch{K}$, $\chuan{h}{\vc}$, $\beta$ and $L$. 
\end{theorem}

We now establish a lower bound in Theorem \ref{th:lower-bound}. 

\begin{theorem}\label{th:lower-bound}
For any $T>0$, $\beta >0$ and $L>0$. Assume that $h \in \mathcal{H}(\beta, L)$, then there exists a constant $C_4>0$ such that for any estimator $\hat{h}_T$ of $h$
\begin{equation}\label{eq:lower-bound}
\underset{h\in\hhoa (\beta, L)}{\sup}\kv \ch{\hat{h}_T - h}^2 \ge C_4 \exp\left(- \frac{2\beta}{2\beta + 1}RT\right).
\end{equation}
\end{theorem}

Contrary to the classical cases of nonparametric estimation (\textit{e.g.} Tsybakov \cite{Tsybakov04}, \ldots), the number of observations $M_T$ is a random variable that converges to $+\vc$ when $T\to +\vc$ which is one of the main difficulty here. From Theorem \ref{th:upper-bound}, when $N_0 > 1$ the upper bound is in $\exp\left(-\frac{2\beta}{2\beta + 1}RT \right)$ which is the same rate as the lower bound. The rate of convergence $\hat h$ is thus optimal. When $N_0 =1$, the upper bound is in $\exp\left(\frac{2\beta}{2\beta+1} \big(-RT + \log(RT)\big)\right)$ that differs with a logarithmic from the rate in the lower bound. The rate of convergence is thus slightly slower than in the case $N_0>1$ and our estimator is optimal up to a logarithmic factor. Furthermore, Theorem 3  illustrates adaptive properties of our procedure: it achieves the rate $\varrho(T)^{-\frac{2\beta}{2\beta+1}}$ over the H\"older class $\hhoa(\beta,L)$ as soon as $\beta$ is smaller than $\beta^*$. So, it automatically adapts to the unknown smoothness of the signal to estimate.

\section{Numerical simulations}
\subsection{Numerical computation of $\hat{h}$}
\indent We use the \textbf{R} software to implement simulations with two original distributions of division kernel $h$ and compare with their estimators. On the interval $[0,1]$, the first distribution to test is $\Beta(2,2)$. $\Beta(a,b)$ distributions on $[0,1]$ are characterized by their densities
\[
h_{\Beta(a,b)}(x) = \frac{x^{a-1}(1-x)^{b-1}}{\mathcal{B}(a,b)}.
\]
where $\mathcal{B}(a,b)$ is the renormalization constant.

\par Since $h$ is symmetric, we only consider the distributions with $a=b$. Generally, asymmetric divisions correspond to $a<1$ and symmetric divisions with kernels concentrated around $\frac 12$ correspond to $a>1$. The smaller the parameter $a$, the more asymmetric the divisions. For the second density, we choose a Beta mixture distribution as
\[
\frac 12\Beta(2,6) + \frac 12\Beta(6,2).
\]
This choice gives us a bimodal density corresponding to very asymmetric divisions.

\begin{figure}[htbp]
\centering
\begin{subfigure}{.47\textwidth}
  \centering
  \includegraphics[scale = 0.35]{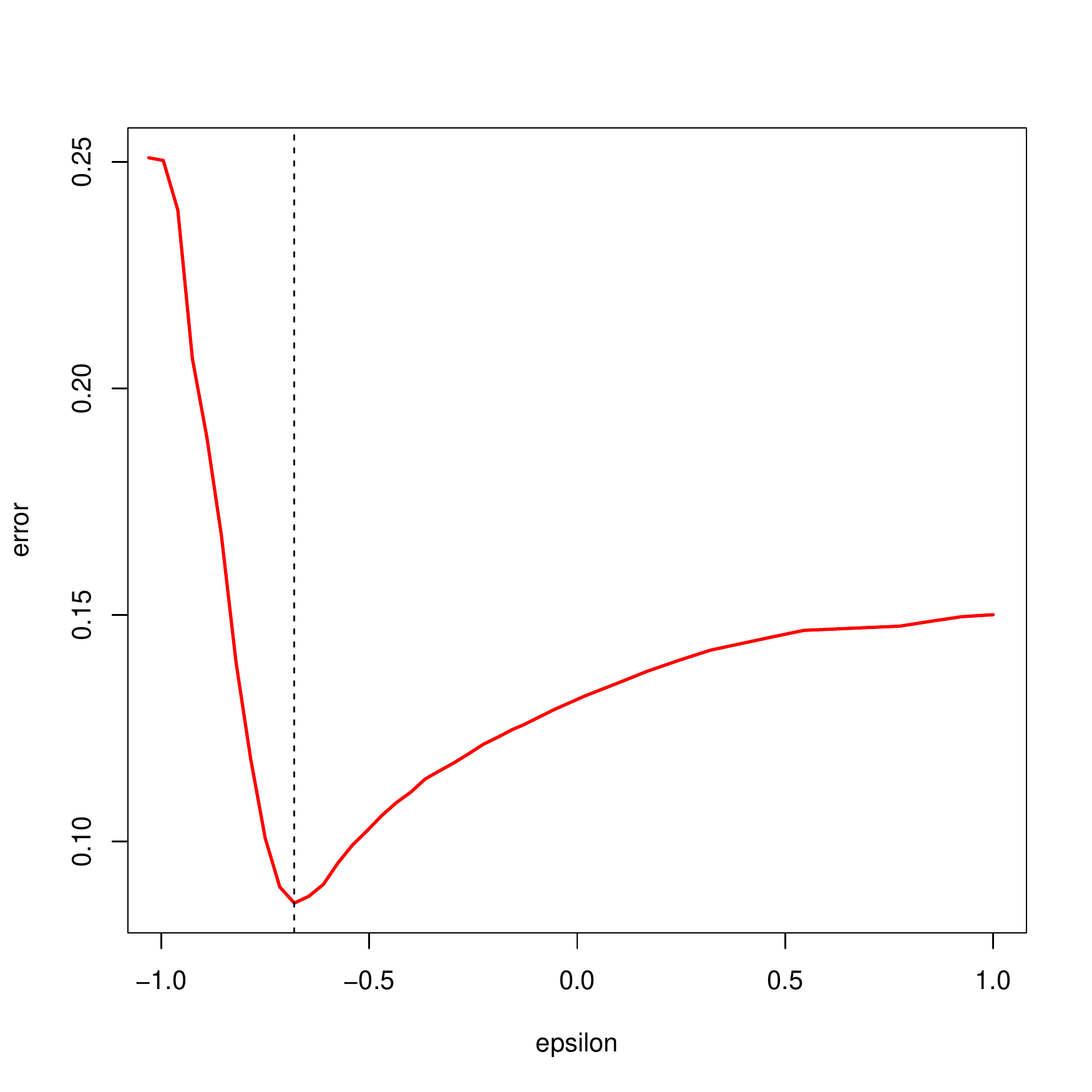}
\caption{}
 \label{fig5b:sub1}
\end{subfigure}%
\begin{subfigure}{.47\textwidth}
  \centering
  \includegraphics[scale = 0.35]{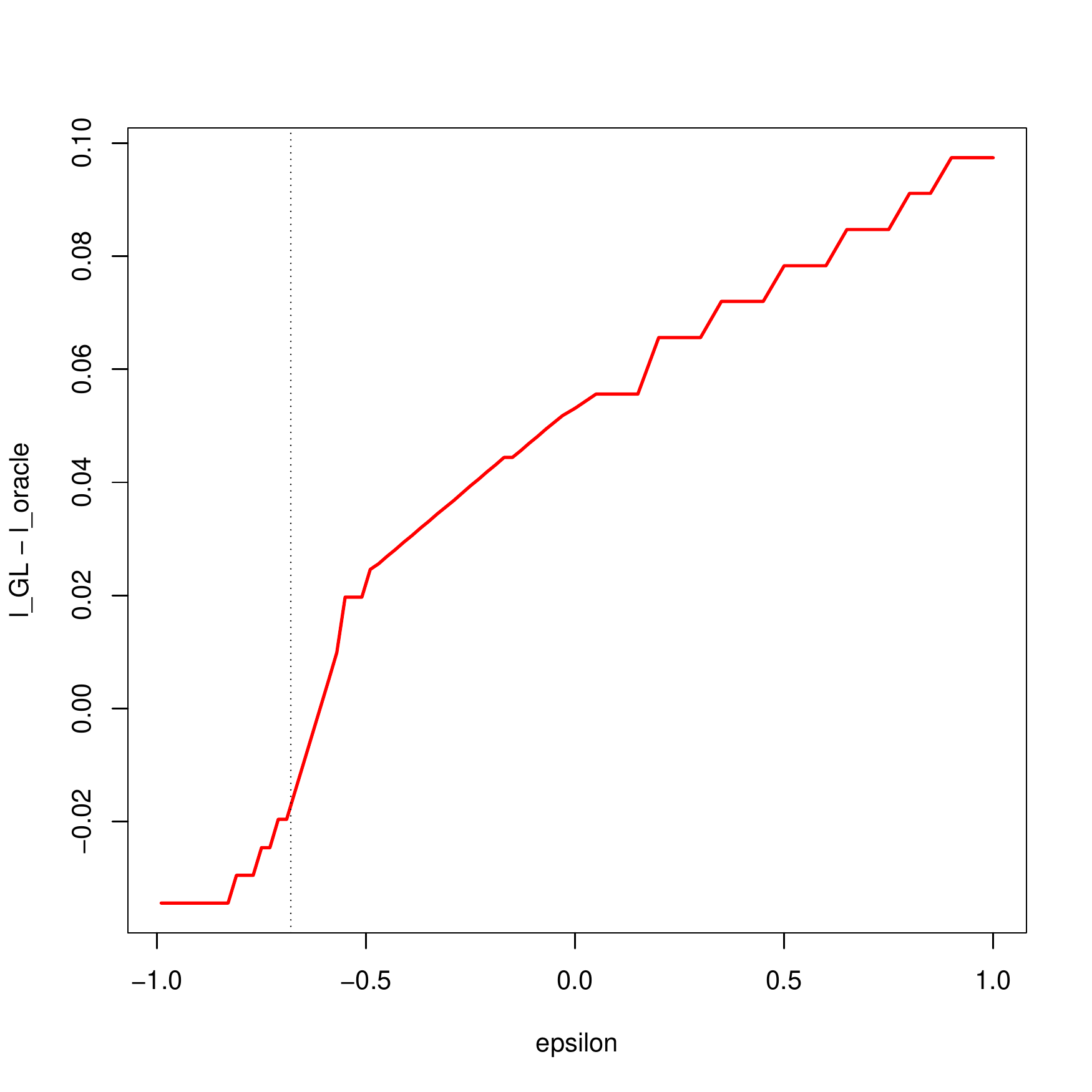}
\caption{}
  \label{fig5b:sub2}
\end{subfigure}
\caption{\itshape (a): MISE's as a function of $\epsilon$. (b): $\hat{\ell} - \ell_{\text{oracle}}$ as a function of $\epsilon$. The dotted lines indicate the optimal value of $\epsilon$ which is used in all simulations. \label{fig5b}}
\end{figure}

\par We estimate $\hat{h}$ by using \eqref{eq:estimator-h} and we take the classical Gaussian kernel $K(x) = (2\pi)^{-1/2}\exp(-x^2/2)$. For the choice of bandwidth, we apply the GL method with the family $H = \left\{1, 2^{-1}, \ldots, \lfloor \delta M_T \rfloor^{-1}\right\}$ for some $\delta>0$ small enough when $M_T$ is large to reduce the time of numerical simulation. We have $\chuan{K}{1}=1$, $\ch{K}=2^{-1/2}\pi^{-1/4}$ and $K_\ell\star K_{\ell'} = K_{\sqrt{\ell^2 + \ell'^2}}$, hence it is not difficult to calculate in practice $\hllp$ as well as $\hlp$. Finally, the value of $\epsilon$ in $\chi = (1+\epsilon)(1 + \chuan{K}{1})$ is chosen to find an optimal value of the MISE. To do this, we implement a preliminary simulation to calibrate $\epsilon$ in which we choose $\epsilon > -1$ to ensure that $1+\epsilon>0$. We compute the MISE and $\hat{\ell} - \ell_{\text{oracle}}$ as functions of $\epsilon$ where $\ell_{oracle} = \argmin_{\ell\in H}\kv\big[\ch{\hat{h}_\ell - h}^2 \big]$ and $h$ is the density of $\Beta(2,2)$. In Figure \ref{fig5b:sub1}, simulation results show that the risk has minimum value at $\epsilon = -0.68$.  This value is not justified from a theoretical point of view. The theoretical choice $\epsilon > 0$ (see Theorem 2) does not give bad results but this choice is too conservative for non-asymptotic practical purposes as often met in the literature (see Bertin et al. \cite{Bertin-Lacour-Rivoirard15} for more details about the GL methodology).  Moreover, following the discussion in Lacour and Massart \cite{Lacour-Massart15} we investigate (see Figure 2b) the difference $\hat{\ell} - \ell_{\text{oracle}}$ and observe some explosions close to $\epsilon=-0.68$. Consequently, we choose $\epsilon = -0.68$ for all following simulations. 

\par
\begin{figure}[htbp]
\centering
\begin{subfigure}{.46\textwidth}
  \centering
  \includegraphics[scale = 0.35]{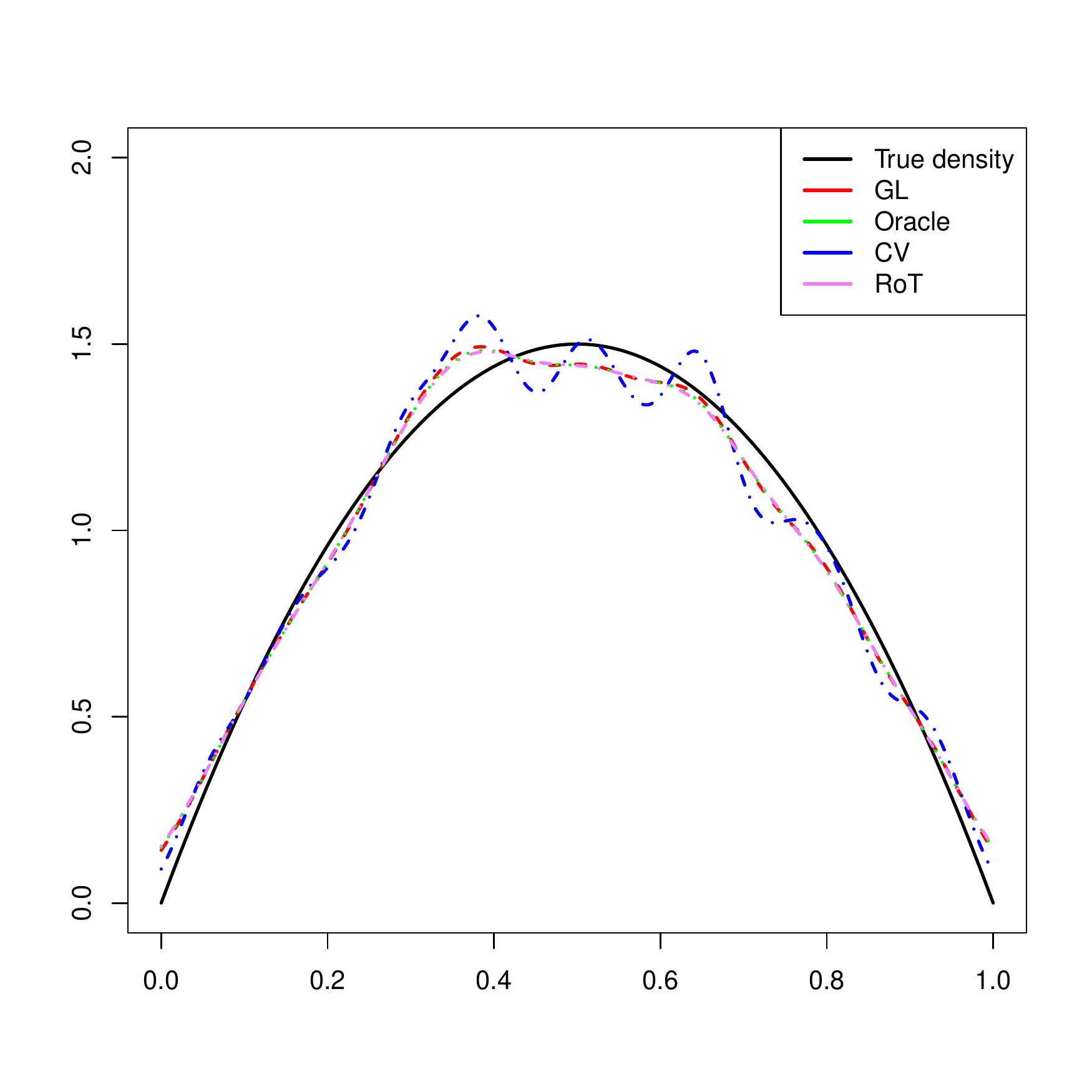}
  \caption{\itshape Reconstruction of $\Beta(2,2)$}
  \label{fig1:sub1}
\end{subfigure}%
\begin{subfigure}{.46\textwidth}
  \centering
  \includegraphics[scale = 0.35]{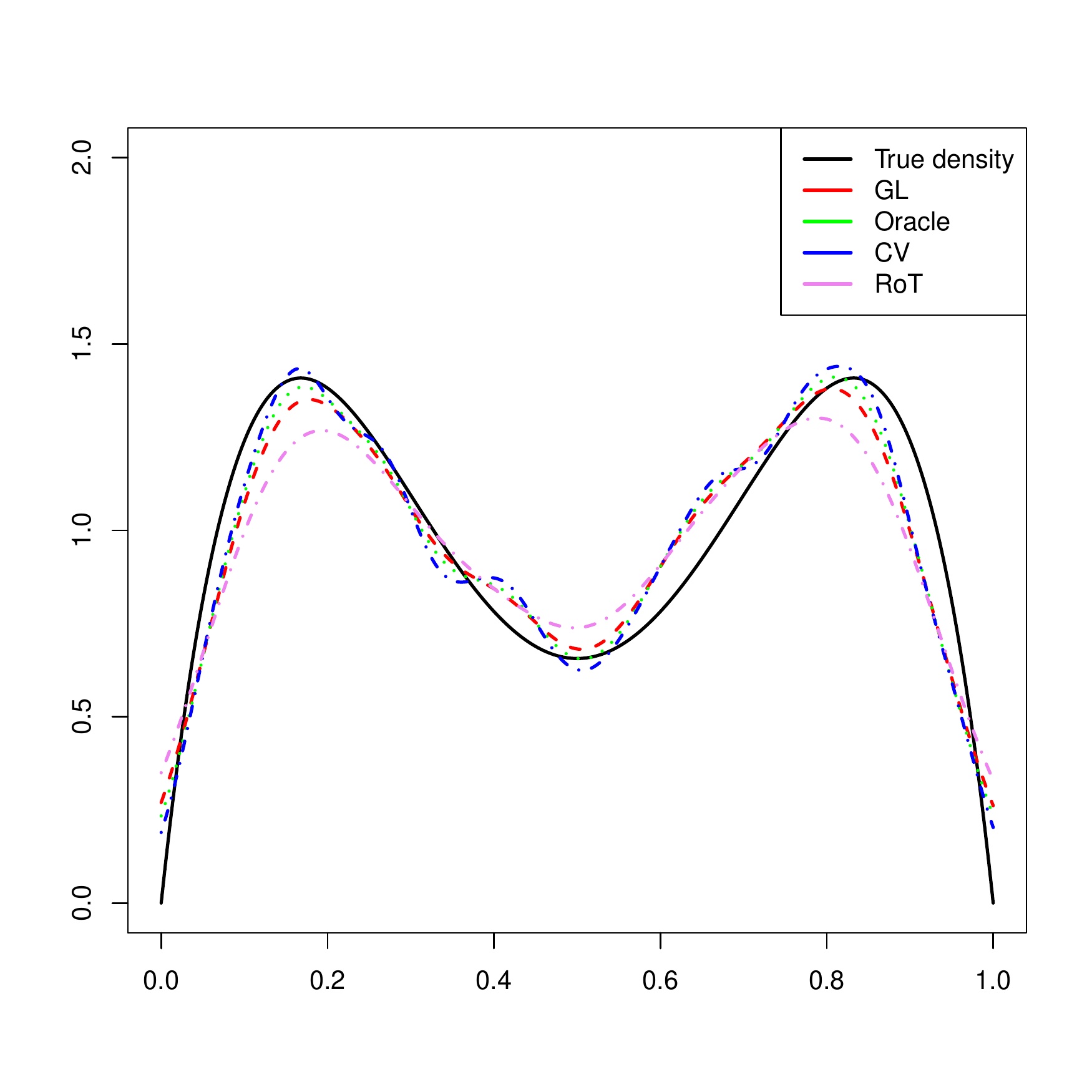}
  \caption{\itshape Reconstruction of beta mixture}
  \label{fig1:sub2}
\end{subfigure}
\caption{\itshape Reconstruction of division kernels with $T=13$.\label{fig1}}
\end{figure}

 Figure \ref{fig1} illustrates a reconstruction for the density of $\Beta(2,2)$ and beta mixture $\frac 12\Beta(2,6) + \frac 12\Beta(6,2)$ when $T=13$. We choose here the division rate and the growth rate $R = 0.5$ and $\alpha = 0.35$ respectively. We compare the estimated densities when using the GL bandwidth with those estimated with the oracle bandwidth. The oracle bandwidth is found by assuming that we know the true density. Moreover, the GL estimators are compared with estimators using the cross-validation (CV) method and the rule of thumb (RoT). The CV bandwidth is defined as follows:
\[
\ell_{CV} = \underset{\ell\in H}{\argmin}\left\{\int \hat{h}_\ell^2(\gamma)d\gamma - \frac{2}{n}\sum_{i=1}^n \hat{h}_{\ell,-i}(\Gamma^1_i) \right\}
\]
where $\hat{h}_{\ell,-i}(\gamma) = \frac{1}{n-1}\sum_{j\neq i}K_\ell\big(\Gamma^1_j - \gamma)$. The RoT bandwidth can be calculated simply by using the formula $\ell_{RoT} = 1.06\hat{\sigma}n^{-1/5}$ where $\hat{\sigma}$ is the standard deviation of the sample $(\Gamma^1_1,\ldots, \Gamma^1_n)$. More details about these methods can be found in Section 3.4 of Silverman \cite{Silverman} or Tsybakov \cite{Tsybakov04}.

\par To estimate the MISE, we implement Monte-Carlo simulations with respect to $T=13, 17$ and $20$. The number of repetitions for each simulation is $\mathcal{M} =100$. Then, we compute the mean of relative error $\bar{e} = (1/\mathcal{M})\sum_{i=1}^\mathcal{M} e_i$ and the standard deviation $\sigma_e = \sqrt{(1/\mathcal{M})\sum_{i=1}^\mathcal{M} (e_i - \bar{e})^2}$ where
\begin{equation}\label{eq:relative-error}
e_i = \frac{\ch{\hat{h}^{(i)} - h}}{\ch{h}},\quad i=1,\ldots,\mathcal{M},
\end{equation} 
and $\hat{h}^{(i)}$ denotes the estimator of $h$ corresponding to $i^{\text{th}}$ repetition.
\begin{table}[ht]
\centering
\begin{tabular}{lllllll}
\hline											
			&				&	GL		&	Oracle	&	CV	&	RoT	& ML method   \\ 
\hline											
$T = 13$	&	$\bar{e}$	&	0.1001	&	0.0840	&	0.1009			&	0.0900	& 0.0610\\
			&	$\sigma_e$	&	0.0585  &	0.0481	&	0.0599			&	0.0577	& 0.0724\\
			&	$\red{\bar{\hat{\ell}}}$	& \red{0.0920}		&\red{0.0845}	&\red{0.0824}	&\red{0.0727} \\

\hline											
$T = 17$	&	$\bar{e}$	&	0.0458	&	0.0397	&	0.0459			&	0.0405	& 0.0166\\
			&	$\sigma_e$	&	0.0260	&	0.0230	&	0.0297			&	0.0237	& 0.0171\\
			&	$\red{\bar{\hat{\ell}}}$	& \red{0.0485}		&\red{0.0497}	&\red{0.0478}	&\red{0.0470} \\
\hline											
$T = 20$	&	$\bar{e}$	&	0.0261	&	0.0241	&	0.0262			& 0.0245	& 0.0088		\\ 
			&	$\sigma_e$	&	0.0140	&	0.0114	&	0.0132			& 0.00121	& 0.0091		\\ 
			&	$\red{\bar{\hat{\ell}}}$	& \red{0.0377}		&\red{0.0359}	&\red{0.0345}	&\red{0.0354} \\
\hline			
\end{tabular}
\caption{\itshape Mean of relative error and its standard deviation for the reconstruction of $\Beta(2,2)$. $\bar{\hat{\ell}}$ is the average of bandwidths for $M=100$ samples. \label{tab1}}
\vspace*{0.1in}
\begin{tabular}{llllll}
\hline											
			&				&	GL		&	Oracle	&	CV	&	RoT	\\
\hline											
$T = 13$	&	$\bar{e}$	&	0.1361	&	0.1245	&	0.1379			&	0.1686	\\
			&	$\sigma_e$	&	0.0672  &	0.0562	&	0.0815			&	0.0537	\\
			&	$\red{\bar{\hat{\ell}}}$	& \red{0.0618}		&\red{0.0527}	&\red{0.0522}	&\red{0.0948} \\

\hline											
$T = 17$	&	$\bar{e}$	&	0.0539	&	0.0534	&	0.0550			&	0.0919	\\
			&	$\sigma_e$	&	0.0180	&	0.0168	&	0.0168			&	0.00223	\\
			&	$\red{\bar{\hat{\ell}}}$	& \red{0.0309}		&\red{0.0272}	&\red{0.0264}	&\red{0.0590} \\
\hline													
\end{tabular}
\caption{\itshape Mean of relative error and its standard deviation for the reconstruction of beta mixture $\frac 12 \Beta(2,6)+\frac 12\Beta(6,2)$.  \label{tab2}}
\end{table}

\par The MISE's are computed for estimated densities using the GL bandwidth, the oracle bandwidth, the CV bandwidth and the RoT bandwidth. For a further comparison, in the reconstruction of $\Beta(2,2)$, we compute the relative error in a parametric setting by comparing the true density $h$ with the density of $\Beta(\hat a,\hat a)$ where $\hat a$ is a Maximum Likelihood (ML) estimator $a$. The simulation results are displayed in Table \ref{tab1} and Table \ref{tab2}. For the density of Beta mixture, we only compute the error with $T=13$ and $T=17$. The boxplot in Figure 4 illustrates the MISE's in Table \ref{tab1} when $T=17$.

\par From Tables \ref{tab1} and \ref{tab2}, we can note that the accuracy of the estimation of $\Beta(2,2)$ and Beta mixture by the GL bandwidth increases for larger $T$. In Figure \ref{fig5}, we illustrate on a log-log scale the mean relative error and the rate of convergence versus time $T$. This shows that the error is close to the exponential rate predicted by the theory. Moreover, we can observe that the errors of Beta mixture are larger than those of $\Beta(2,2)$ with the same $T$ due to the complexity of its density. In both cases, the error estimated by using oracle bandwidth is always smaller. The GL error is slightly smaller than the CV error. The RoT error can show very good behavior but lacks of stability.  Overall, we conclude that the GL method has a good behavior when compared to the cross validation method and rule-of-thumb. As usual, we also see that the ML errors are quite smaller than those of nonparametric approach but the magnitude of the mean $\bar e$ remains similar.\\

\begin{figure}[h]
\centering
\includegraphics[scale = 0.32]{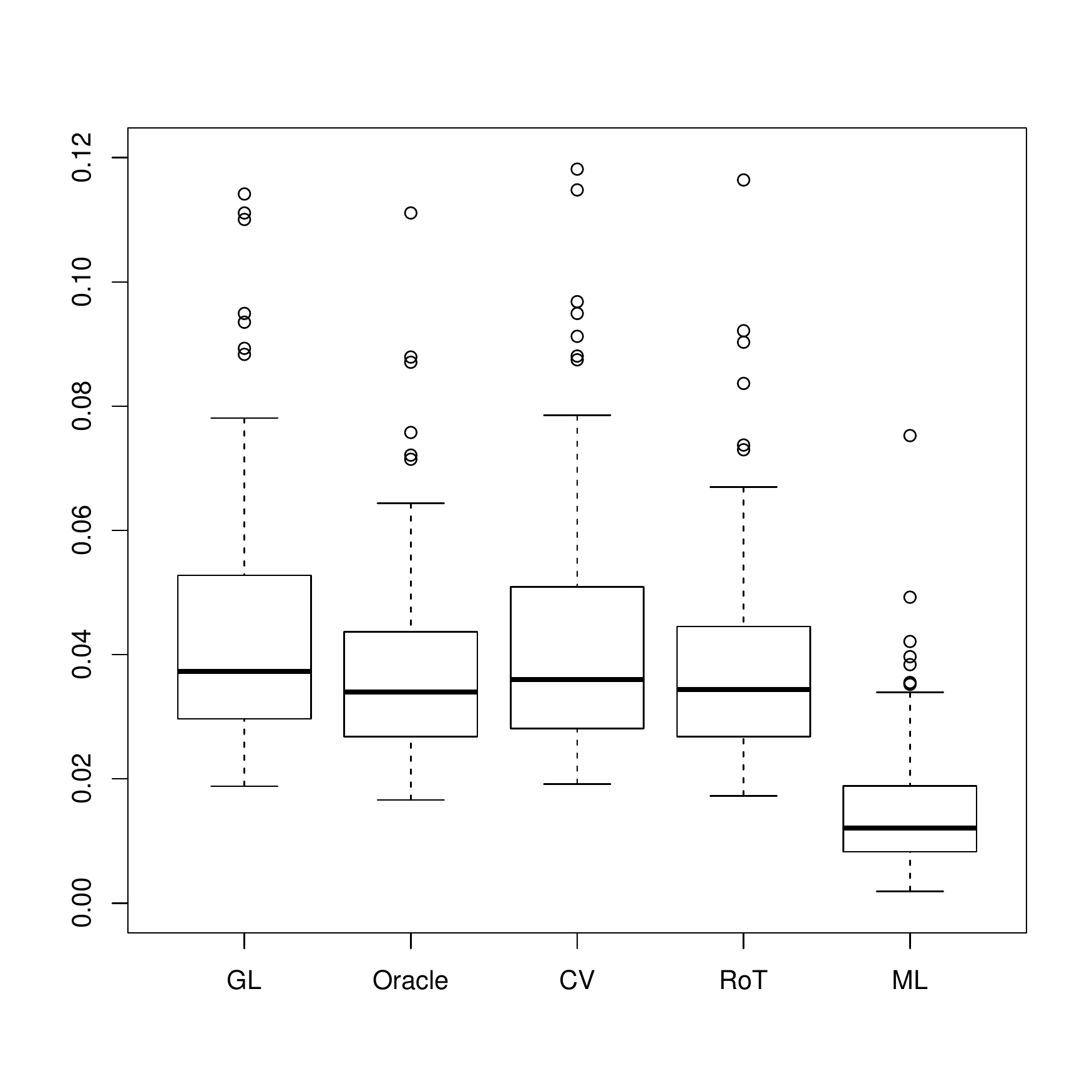}
\vspace{-1.5em}
\caption{\itshape Errors of estimated densities of $\Beta(2,2)$ when $T=17$. \label{fig4}}
\end{figure}
\begin{figure}[htbp]
\centering
\includegraphics[scale = 0.32]{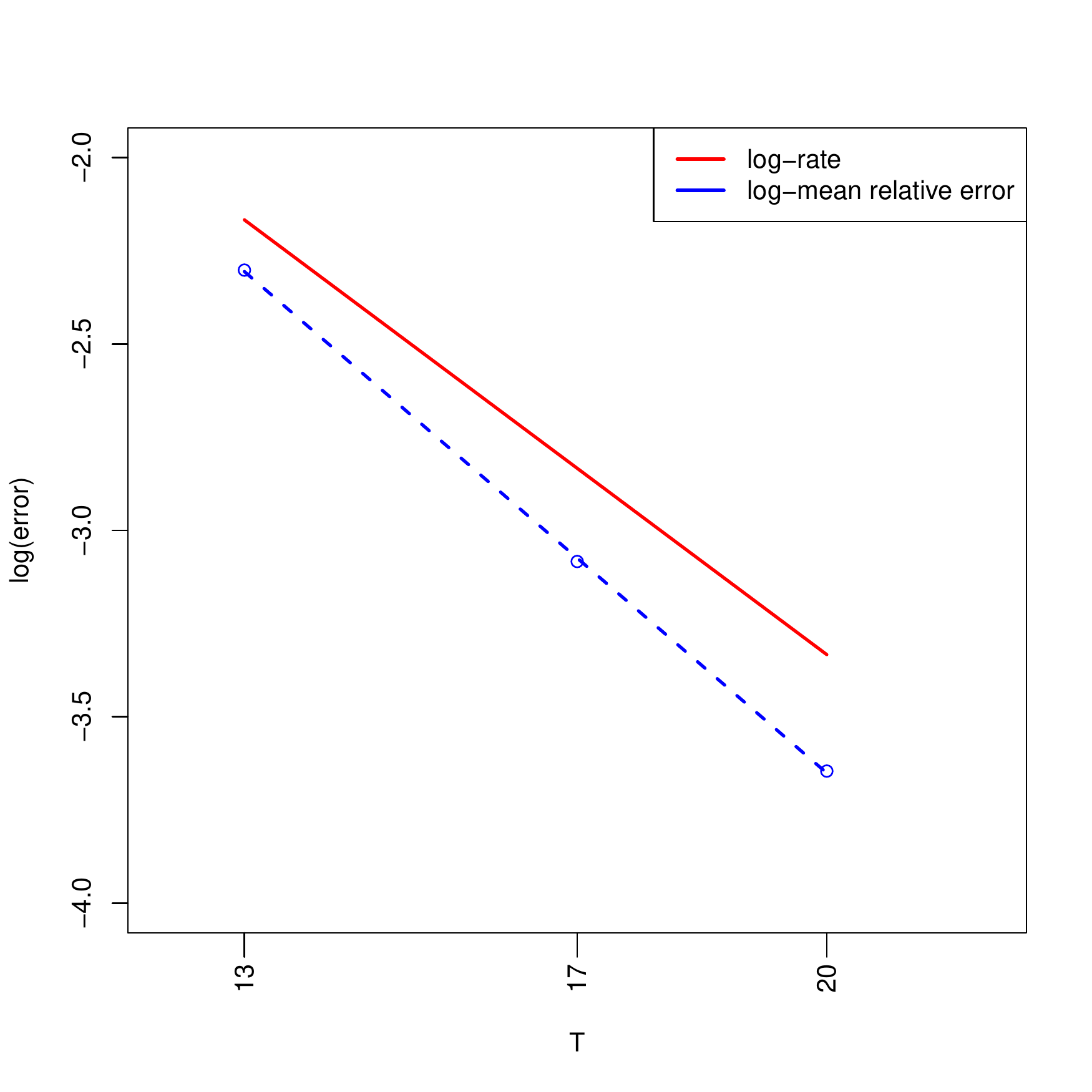}
\vspace{-1.5em}
\caption{\itshape The log-mean relative error for the reconstruction of $\Beta(2,2)$ compared to the log-rate (solid line) computed with $\beta=1$. \label{fig5}}
\end{figure}
\newpage
Since $h$ is symmetric on $[0,1]$ with respect to $\frac 12$, the estimator $\hat{h}$ can be improved and we can introduce
\[
\tilde{h}(x) = \frac 12\left(\hat{h}(x) + \hat{h}(1-x) \right),
\]
which is symmetric by construction and satisfies also \eqref{eq:upper-bound}. We compute the mean of relative error for the estimator $\tilde{h}$ with the estimation of $\Beta(2,2)$ and Beta mixture. The results are displayed in Table \ref{tab3}. Compared with the error in Table \ref{tab1} and \ref{tab2}, one can see as expected that the errors for the reconstruction of $\tilde{h}$ are smaller. However, these errors are of the same order, indicating that the estimator $\hat{h}$ had already good symmetric properties.

\begin{table}[htbp]
\centering
\begin{tabular}{llllll}
\hline											
				&			&	GL		&	Oracle	&	CV				&	RoT	\\ 
\hline											
$\Beta(2,2)$	&	$T=13$	&	0.0785	&	0.0634	&	0.0762			&	0.0644	\\
				&	$T=17$	&	0.0356  &	0.0309	&	0.0356			&	0.0309	\\
\hline											
Beta mixture	&	$T=13$	&	0.1117	&	0.0953	&	0.1030			&	0.1584	\\
				&	$T=17$	&	0.0450	&	0.0414	&	0.0417			&	0.0893	\\
\hline											
\end{tabular}
\caption{\itshape Mean of relative error for the reconstruction of $\tilde{h}$. \label{tab3}}
\end{table}

\subsection{Influence of the distribution on the mean age}
\begin{figure}[htbp]
\centering
\includegraphics[scale = 0.4]{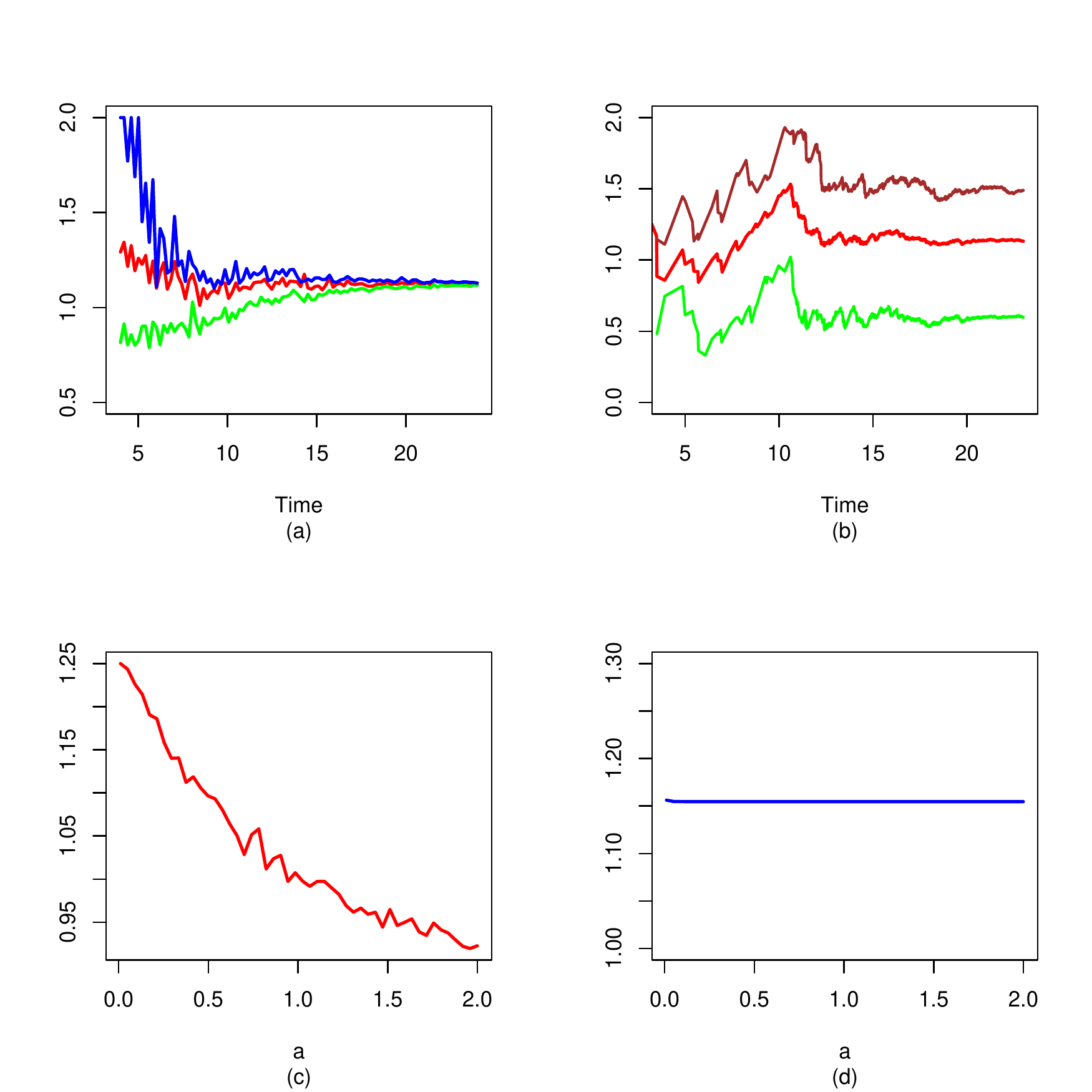}
\caption{\itshape (a) Average mean, $1^{\text{st}}$ and $3^{\text{rd}}$ quartiles for the sample of means for $50$ trees. (b) Average mean, $1^{\text{st}}$ and 
$3^{\text{rd}}$ quartiles for one tree. (c) Average of $Q_{75} - Q_{25}$ with $a\in [0,2]$ at $t=12$. (d) Mean age with $a\in [0,2]$ at $t=12$.  \label{fig3}}
\end{figure}

\noindent For $t\ge 0$, recall the mean age defined in \eqref{eq:mean-age}. To study the influence of the distribution on the mean age, we simulate $n=50$ trees with respect to $t=6, 6+\triangle t,\ldots, 24$ with $\triangle t = 0.36$. For each sample ($\bar x^{(1)}_t$, \ldots, $\bar x^{(n)}_t$), we compute the average mean, the $1^{\text{st}}$ ($Q_{25}$) quartile and $3^{\text{rd}}$ ($Q_{75}$) quartile. Figure \ref{fig3}a and \ref{fig3}b show the simulation results corresponding to the density of $\Beta(2,2)$ with $\alpha = 0.45$ and $R =0.4$. One can see that the average of mean age and the mean age converge to $\frac{\alpha}{R} = 1.125$ for larger $t$. This agrees with the theoretical result proved in Section \ref{sec:mean-age}. 

\par Moreover, $Q_{25}$ and $Q_{75}$ vary when the parameter $a$ changes. In Figure \ref{fig3}c, we draw a fitted curve of the average of 
$(Q_{75} - Q_{25})$ when $a$ varies from $0$ to $2$. As we mentioned in the introduction, if divisions are more asymmetric corresponding to small values of $a$, the toxicities concentrate on few cells, \textit{i.e.} we have more older cells after the divisions. This explains the decreasing trend in the average of $(Q_{75} - Q_{25})$. Finally, Figure \ref{fig3}d displays the average of mean ages with respect to $a$. One can note that it does not change when we replace the kernel distribution, \textit{e.g} $\Beta(0.6, 0.6)$ instead of $\Beta(2,2)$.

\section{Proofs}

\subsection{Proof of Proposition \ref{prop:N_T}}
\textit{ii)} The proof of ii) can be found easily in literature. Here we refer to \cite{Ross95}, Section 5.3 for this proof. \\[6pt]
\textit{i)} Let us prove that $\lim_{T\to +\vc} N_T = \lim_{j\to+\vc} N_{T_j} = +\vc$. Since our model has only births and no death, 
$(N_t)_{t\in [0,T]}$ is a non-decreasing process: $N_{T_j}=N_0 + j$. All the $T_j$'s are finite and $\lim_{j\to +\vc} N_{T_j} = +\vc$ a.s. From \textit{ii)}, we have $\kv[N_T] = N_0 e^{RT}$. Hence, we deduce from the estimate $\usup{t\in [0,T]}\kv[N_t]<+\vc$ for all $T>0$ that $T_j 
\underset{j\to+\vc}{\longrightarrow}+\vc$ a.s. Then we also have $\lim_{T\to +\vc} N_T = +\vc$ a.s.\\[6pt]
\textit{iii)} Let $p= e^{-RT}$. When $N_0 = 1$, $N_T\sim \Geom(p)$. Then we have
\begin{align*}
\kv\left[\frac{1}{N_T}\right] &= \sum_{n=1}^{\infty} \frac{1}{n}\xs\big(N_T = n\big) = \sum_{n=1}^{\infty} \frac{1}{n} p(1-p)^{n-1} \\
&= \frac{p}{1-p}\sum_{n=1}^{+\vc} \frac{(1-p)^n}{n} = -\frac{p}{1-p}\log(p).
\end{align*}
Replace $p$ with $e^{-RT}$, we obtain \eqref{eq:inverse-N_T-geom}.\\
When $N_0 >1$, $N_T\sim \mathcal{NB}(N_0, p)$. Hence, we have
\begin{align}
\kv\left[\frac{1}{N_T}\right] =& \sum_{n = N_0}^{\infty} \frac{1}{n}\binom{n-1}{n - N_0} p^{N_0}( 1- p)^{n-N_0}\notag\\
=& \left(\frac{p}{1-p}\right)^{N_0}\sum_{n = N_0}^{\infty} \frac{1}{n}\binom{n-1}{n - N_0} (1-p)^n \notag\\
:=& \left(\frac{p}{1-p}\right)^{N_0}f(1-p), \label{eq:temp4}
\end{align}
where $f(x)=\sum_{n=N_0}^{+\infty} \frac{1}{n} { n-1 \choose n-N_0}  x^n$. We can differentiate $f(x)$ by taking derivative under the sum. Then:
\begin{align*}
\frac{d}{dp}f(1-p) = &  - \sum_{n=N_0}^{+\infty}  { n-1 \choose n-N_0}  (1-p)^{n-1}\\
= &  - \frac{(1-p)^{N_0-1}}{p^{N_0}} \sum_{n=N_0}^{+\infty}  { n-1 \choose n-N_0}  p^{N_0}(1-p)^{n-N_0}=  - \frac{1}{p}\Big(\frac{1}{p}-1\Big)^{N_0-1},
\end{align*}since the sum is 1 (we recognize the negative binomial).\\
Hence,
\begin{align}
\frac{d}{dp}f(1-p) =& -\frac{1}{p}\left[\sum_{k=1}^{N_0-1} \binom{N_0-1}{k}\frac{1}{p^k}(-1)^{N_0-1-k} + (-1)^{N_0-1}\right] \notag \\
=& (-1)^{N_0}\left[\sum_{k=1}^{N_0-1}\binom{N_0-1}{k} \frac{(-1)^k}{p^{k+1}} + \frac 1p \right].\label{eq:temp-1}
\end{align}
Integrating equation \eqref{eq:temp-1} and notice that $f(0) = 0$, we get
\begin{align}
f(1-p) = & (-1)^{N_0}\left[\sum_{k=1}^{N_0-1} \binom{ N_0 - 1}{k } \frac{(-1)^k}{k}\Big(-\frac{1}{p^k}\Big) + \log(p) \right] \notag\\
=& (-1)^{N_0-1}\left[\sum_{k=1}^{N_0-1} \binom{ N_0 - 1}{k } \frac{(-1)^k}{k}\frac{1}{p^k} + \log\left(\frac 1p \right) \right]. \label{eq:temp5}
\end{align}
Combine \eqref{eq:temp4},\eqref{eq:temp5} and replace $p$ with $e^{-RT}$, we get \eqref{eq:inverse-N_T-nb}.

\noindent\textit{iv)} We first prove the lower bound of \eqref{eq:bound-inverse-N_T-nb}. From \eqref{eq:SDE-model}, taking $f_t(x) = 1$, we have 
\begin{equation}\label{eq:N_T}
 N_T = N_0 + \int_0^T\int_{\ehoa} \id_{\{i\le N_{s-}\}}Q(ds,di,d\gamma).
\end{equation}
Applying It\^{o} formula for jump processes (see \cite{Ikeda89}, Theorem 5.1 on p.67) to \eqref{eq:N_T}, we obtain 
\begin{align*}
 \frac{1}{N_T} &= \frac{1}{N_0} + \int_0^T\int_{\ehoa}\left(\frac{1}{N_{s-}+1} -  \frac{1}{N_{s-}}\right)\id_{\{i\le N_{s-}\}}Q(ds,di,d\gamma)\\
&= \frac{1}{N_0} - \int_0^T\int_{\ehoa} \frac{1}{N_{s-}\left(N_{s-} + 1 \right)}\id_{\{i\le N_{s-}\}}Q(ds,di,d\gamma).
\end{align*}
Hence,
\begin{equation}\label{eq:temp0}
\kv\left[\frac{1}{N_T} \right] = \frac{1}{N_0} - \kv\left[\int_0^T\frac{1}{N_{s}\left(N_{s} + 1 \right)}RN_s ds \right]
= \frac{1}{N_0} - R\int_0^T\kv\left[\frac{1}{N_s + 1} \right]ds.
\end{equation}
Since $N_s \ge N_0$, we have
$
\frac{1}{N_s + 1} \le \frac{1}{N_s} .
$
Therefore, (\ref{eq:temp0}) implies that
\begin{equation}\label{eq:temp1}
\kv\left[\frac{1}{N_T} \right] \ge \frac{1}{N_0} - R\int_0^T\kv\left[\frac{1}{N_s} \right]ds  .
\end{equation}
By comparison of $\kv\left[\frac{1}{N_T}\right]$ with the solutions of the ODE
$
\frac{d}{dT}u(T) = - Ru(T)
$
with $u(0) = 1/N_0$, we finally obtain
\[
\kv\left[\frac{1}{N_T} \right] \ge \frac{1}{N_0}e^{-RT}.
\] 

For the upper bound, notice that $\kv\big[\frac{1}{N_T}\big] \le \kv\left[\frac{1}{N_T - 1} \right]$ for $N_0 > 1$. Then we have
\begin{align*}
\kv\Big[\frac{1}{N_T-1}\Big]= & \sum_{n=N_0}^{+\infty} \frac{1}{n-1} \binom{n-1}{n-N_0} p^{N_0} (1-p)^{n-N_0}\\
 = & \sum_{n=N_0}^{+\infty} \frac{(n-2)!}{(n-N_0)!(N_0-1)!} p^{N_0} (1-p)^{n-N_0}\\
 = & \frac{p}{N_0-1} \sum_{n=N_0}^{+\infty} \frac{(n-2)!}{(n-N_0)!(N_0-2)!} p^{N_0-1} (1-p)^{n-N_0}\\
 = & \frac{p}{N_0-1} \sum_{m=N_0-1}^{+\infty} \frac{(m-1)!}{(m-(N_0-1))!((N_0-1)-1)!} p^{N_0-1} (1-p)^{m-(N_0-1)}\\
 = & \frac{p}{N_0-1} = \frac{e^{-RT}}{N_0 - 1},
\end{align*}
by changing the index in the sum ($m=n-1$) and by recognizing the negative binomial with parameter $(N_0-1,p)$.
Hence, we conclude that for $N_0>1$
\[
\frac{e^{-RT}}{N_0} \le \kv\left[ \frac{1}{N_T}\right] \le \frac{e^{-RT}}{N_0-1}.
\]
This ends the proof of Proposition \ref{prop:N_T}.

\subsection{Proof of Lemma \ref{lem1}}\label{proof:lem-mean-age}

By symmetry of $h$ with respect to $1/2$, we have:
\begin{align}
Y_t &= Y_0 + \int_0^t\left(\alpha + 2R\int_0^1\left(\gamma Y_s - Y_s \right)h(\gamma)d\gamma \right)ds + U_t \notag \\
&= Y_0 + \int_0^t\left(\alpha - 2R Y_s\int_0^1\gamma h(\gamma)d\gamma \right)ds + U_t \notag \\
&= Y_0 + \int_0^t\left(\alpha - R Y_s\right)ds + U_t. \notag
\end{align}
where $U_t$ is a square-integrable martingale. \\
Let $\tilde{Y}_t = Y_t e^{Rt}$, $\tilde Y_0 = Y_0$. By It\^o formula, we get
\[
\tilde{Y}_t = \tilde{Y_0} + \frac{\alpha}{R}\left(e^{Rt} - 1\right) + \int_0^t e^{Rs}dU_s.
\]
Replacing $\tilde{Y}_t$ by $Y_te^{Rt}$, we obtain 
\[
Y_t = \left(Y_0 - \frac{\alpha}{R} \right)e^{-Rt} + \frac{\alpha}{R} + \int_0^t e^{-R(t-s)}dU_s.
\]
We end the proof by taking the expectation and the limit as $t\rightarrow +\vc$ of $Y_t$ to obtain \eqref{eq:ave-mean-age} and \eqref{eq:limit-ave-mean-age}. 

\subsection{Proof of Theorem \ref{th:mean-age}}\label{proof:th-mean-age}
We will show that the process $Y$ satisfies ergodicity and integrability assumptions in Bansaye {et al.} \cite{Bansaye11a} (see (H1) - (H4), Section 4). More precisely:
\begin{itemize}
\item[1.] $\kv\left[Y_t\right] < +\vc$ for all $t\ge 0$.
\item[2.] There exists $\varpi < R$  and $c>0$ such that $\kv\left[Y_t^2 \right] < ce^{\varpi t}$ for all $t\ge 0$.
\end{itemize}

From \eqref{eq:ave-mean-age} we note that $\kv[Y_t] <+\vc$ for all $t\ge 0$. To prove the second point, from \eqref{eq:infinitesimal-generator} we have
\begin{align}
\kv[Y_t^2] &= \kv\left[Y_0^2 + \int_0^t\left(2\alpha Y_s +2R\int_0^1\left(\gamma^2Y_s^2 - Y_s^2 \right)h(\gamma)d\gamma \right)ds \right]\notag\\
&=Y_0^2 + 2\alpha\int_0^t\kv[Y_s]ds - 2\theta R\int_0^t\kv[Y_s^2]ds,\label{eq:mean-age-temp3}
\end{align}
with $\theta = \int_0^1(1-\gamma^2)h(\gamma)d\gamma$ and $0<\theta <1$.

\par Substituting $\kv[Y_t] = (Y_0 - \alpha/R)e^{-Rt} + \alpha/R$ into \eqref{eq:mean-age-temp3}, we see that $\kv(Y_t^2)$ solves the following equation:
\begin{equation}\label{eq:mean-age-temp4}
\frac{d\kv[Y_t^2]}{dt} = -2\theta R\kv[Y_t^2] + \left(2\alpha Y_0 - \frac{2\alpha^2}{R}\right)e^{-Rt} + \frac{2\alpha^2}{R}.
\end{equation}
The solution of the equation \eqref{eq:mean-age-temp4} is:
\begin{equation}\label{eq:mean-age-temp5}
\kv[Y_t^2] = e^{-2\theta Rt}\left[Y_0^2 + \int_0^t e^{2\theta Rs}\left(\big(2\alpha Y_0 - \frac{2\alpha^2}{R}\big)e^{-Rs} + \frac{2\alpha^2}{R} \right)ds \right].
\end{equation}
Hence, if $\theta = \frac 12$, we have
\begin{align*}
\kv[Y_t^2] &= Y_0^2e^{-Rt} + \big(2\alpha Y_0 - \frac{2\alpha^2}{R} \big)te^{-Rt} + \frac{2\alpha^2}{R^2}\left(1-e^{-Rt}\right) \\
&\le Y_0^2e^{-Rt} + \big(2\alpha Y_0 - \frac{2\alpha^2}{R} \big)e^{-(R-\theta)t} + \frac{2\alpha^2}{R^2} \\
&\le \left(Y_0^2 + 2\alpha Y_0 + \frac{2\alpha^2}{R} + \frac{2\alpha^2}{R^2}\right)e^{(0\vee (\theta-R))t} = c_1 e^{\varpi t},
\end{align*}
with $\varpi = 0\vee (\theta-R):= \max(0,\theta-R)$.\\[6pt]
If $\theta\neq \frac 12$, 
\begin{align}
\kv[Y_t^2] &= e^{-2\theta Rt}\left[Y_0^2 + 
\big(2\alpha Y_0 - \frac{2\alpha^2}{R}\big)\int_0^t e^{(2\theta-1)Rs}ds + 
\frac{2\alpha^2}{R}\int_0^t e^{2\theta Rs}ds \right] \notag \\
&= Y_0^2e^{-2\theta Rt} + \big(2\alpha Y_0 - \frac{2\alpha^2}{R}\big)\frac{1}{(2\theta-1)R}\left(e^{-Rt}-e^{-2\theta Rt}\right)\notag\\
&\phantom{aaaaaaaaaaaaaaaaaaaaaaaaaaaaaaaaaaaa}  + \frac{\alpha^2}{\theta R^2}\left(1-e^{-2\theta Rt}\right)\notag \\
&\le \left(Y_0^2 + \big(2\alpha Y_0 + \frac{2\alpha^2}{R}\big)\frac{1}{|2\theta-1|R} + \frac{\alpha^2}{\theta R^2} \right) = c_2. \notag
\end{align}
Thus, if we set $c = \max(c_1,c_2)$ then $\kv\left[Y_t^2 \right] < ce^{\varpi t}$ for all $t\ge 0$.\\[6pt]
The infinitesimal generator $A$ of $Y$ is defined for $\mathcal{C}^1$ test functions as
\[
 Af(x) = \alpha f'(x) + 2R\int_0^1 \left(f(\gamma x) - f(x) \right)h(\gamma) d\gamma.
\]
For $V(x)=x$ and $f(x)=x+1$, we have 
$$AV(x)=\alpha - Rx  \leq -\frac{R}{2} f(x)+\Big(\alpha +\frac{R}{2}\Big) \id_{\big\{x\leq \frac{2\alpha}{R}+1\big\}}.$$
Hence, by Theorem 5.3 of Meyn and Tweedie \cite{Meyn_Tweedie93}, there exists $\pi\in\mhoa_F(\rhoad)$ such that $\lim_{t\to+\vc} E[Y_t] = \langle \pi, f \rangle = \frac{\alpha}{R}$. Finally, applying Theorem 4.2 of \cite{Bansaye11a}, we obtain the result
\[
\underset{t\to +\vc}{\lim} \frac{\langle Z_t, f\rangle}{N_t} = \langle \pi, f\rangle = \frac{\alpha}{R}.
\]

\subsection{Proof of Proposition \ref{prop:mean-var-h}}\label{proof:prop-mean-var-h}

To prove \eqref{eq:mean-h}, let us remark that the number of random divisions $M_T$ is independent of $(\Gamma_i^1)_{i\in\nhoa}$, because the division rate $R$ is constant and because of the construction of our stochastic process. Therefore, we have
\begin{align*}
\kv\big[\hl |M_T \big] &= \kv\Big[\frac{1}{M_T}\sum_{i=1}^{M_T} \Kl \big| M_T\Big]= \frac{M_T\kv[\Kll]}{M_T} \\
&= \kv\big[\Kll \big] = K_\ell\star h (\gamma),
\end{align*}
and $\kv\big[\hl \big] = \kv\left[\kv\big[\hl |M_T \big] \right] = K_\ell\star h (\gamma)$. By similar calculations as \eqref{eq:mean-h}, we obtain \eqref{eq:var-cond-h} and \eqref{eq:var-h}.

\par To prove \textit{ii)}, by the Strong Law of Large Numbers, we have
\[
\frac 1n \sum_{i=1}^n K_\ell(\gamma - \Gamma^1_i) \overset{\text{a.s}}{\longrightarrow} \kv\big[K_\ell(\gamma - \Gamma^1_1 )\big] 
\quad\text{ as } n\to +\vc. 
\]

From \eqref{eq:limit-N_T}, we have $\lim_{T\to+\vc} N_T = +\vc$ (a.s). Since $M_T = N_T - N_0$ and $N_0$ is deterministic, this yields
\[
 \frac{1}{M_T}\sum_{i=1}^{M_T} K_\ell\big(\gamma - 
\Gamma^1_i\big) \overset{\text{a.s}}{\longrightarrow} \kv\big[K_\ell(\gamma - \Gamma^1_1) \big] = K_\ell\star h(\gamma).
\]

This ends the proof of Proposition \ref{prop:mean-var-h}.

\subsection{Proof of Proposition \ref{prop:bias-variance-decomposition}}\label{proof:prop-bias-variance-decomposition}

We have
\begin{align*}
\kv\big[\ch{\htl - h}|M_T \big] \le \ch{h - K_\ell\star h} + \kv\big[\ch{\htl - \kv[\htl]}|M_T \big].
\end{align*}
For the variance term, using that $\kv\big[\hat{h}_{\ell}(\gamma) \big] = \kv\big[\hat{h}_{\ell}(\gamma)|M_T \big]$ 
\begin{align}
\kv\big[\ch{\htl - \kv[\htl]}^2 |M_T\big] &= 
\kv\Big[\int_\rhoa \big|\htl(\gamma) - \kv\big[\htl(\gamma)\big]  \big|^2d\gamma\big|M_T \Big] \notag \\
&= \int_\rhoa \ps\Big[\htl(\gamma) \big| M_T \Big]d\gamma \notag \\
&= \frac{1}{M_T}\int_\rhoa \ps\Big[K_\ell(\gamma - \Gamma^1_1) \Big]d\gamma \notag\\
&\le \frac{1}{M_T}\int_\rhoa \kv\Big[K_\ell^2(\gamma - \Gamma^1_1)\Big]d\gamma \notag
\end{align}
By Fubini's theorem, we get
\begin{align*}
\int_\rhoa \kv\Big[K_\ell^2(\gamma - \Gamma^1_1)\Big]d\gamma &= \int_\rhoa\int_\rhoa K^2_\ell(\gamma - u)h(u)du\,d\gamma \\
&= \int_\rhoa h(u)\left(\int_\rhoa K^2_\ell(\gamma - u)d\gamma\right)du \\
&= \ch{K_\ell}^2\int_\rhoa h(u)du = \frac{\ch{K}^2}{\ell}.
\end{align*}
Then we have
\begin{equation}
\kv\big[\ch{\htl - \kv[\htl]}^2 |M_T\big] \le  \frac{\ch{K}^2}{M_T\ell}. \label{eq:upper-bound-var}
\end{equation}
Hence, applying Cauchy-Schwarz's inequality, we obtain \eqref{eq:bias-variance-decomposition}.
This ends the proof of Proposition  \ref{prop:bias-variance-decomposition}.

\subsection{Proof of Theorem \ref{th1}}\label{proof-oracle-inequality}

This proof is inspired by the proof of Doumic {et al.} \cite{Doumic12}. However, our problem here is that the number of observations $M_T$ is random.
To overcome this difficulty, we work conditionally to $M_T$ to get concentration inequalities. \\[3pt] 
\indent Hereafter, we refer $\int f$ to $\int_\rhoa f$ and since the support of 
$h$ is $(0,1)$, we can write $\int h(\gamma)d\gamma$ instead of $\int_0^1 h(\gamma)d\gamma$. Recall that
\[
A(\lb):=\underset{\lb'\in H}{\sup}\Big\{\ch{\hat h_{\lb,\lb'} - \hat h_{\lb'}} 
- \frac{\chi\ch{K}}{\sqrt{M_T\lb'}}\Big\}_{+}.
\]

\noindent Then, for any $\ell \in H$, we have
\[
\ch{\hat h - h} \le A_1 + A_2 + A_3,
\]
where
\begin{align*}
A_1 &:= \ch{\hat{h} - \hlml} \le A(\ell) + \frac{\chi\ch{K}}{\sqrt{M_T\hat{\ell}}}, \\
A_2 &:= \ch{\hlml - \hl} \le A(\hat{\ell}) + \frac{\chi\ch{K}}{\sqrt{M_T \ell}}, \\
A_3 &:= \ch{\hl - h}.
\end{align*}
By definition of $\hat\ell$, we have
\begin{equation}\label{eq:A1A2}
A_1 + A_2 \le 2A(\ell) + 2\frac{\chi\ch{K}}{\sqrt{M_T \ell}},
\end{equation}
and
\begin{align}\label{Aell}
A(\ell) &\le \underset{\ell'\in H}{\sup}\left\{
	\ch{\big(\hllp - \kv[\hllp] \big) - \big(\hlp - \kv[\hlp] \big)} \right. \notag\\
&\phantom{aaaaaaaaaaaaaaaaa}\left. + \ch{\kv[\hllp] - \kv[\hlp]} - \frac{\chi\ch{K}}{\sqrt{M_T \ell'}}\right\}_{+}\notag\\
&\le \xi_T(\ell) + \underset{\ell'\in H}{\sup}\left\{\ch{\kv[\hllp] - \kv[\hlp]}\right\},
\end{align}
where
\begin{equation}\label{eq:proof-th1-xi}
\xi_T(\ell) = \underset{\ell'\in H}{\sup}\left\{
	\ch{\big(\hllp - \kv[\hllp] \big) - \big(\hlp - \kv[\hlp] \big)}
	- \frac{\chi\ch{K}}{\sqrt{M_T\ell'}} \right\}_{+}.
\end{equation}

\noindent For the term $\underset{\ell'\in H}{\sup}\left\{\ch{\kv[\hllp] - \kv[\hlp]}\right\}$, we have
\begin{align}
\kv[\hllp] &- \kv[\hlp] = \int \big(K_\ell\star K_{\ell'}\big)(\gamma - u)h(u)du 
	- \int K_{\ell'}(\gamma - v)h(v)dv \notag \\
&= \int\int K_\ell(\gamma - u - t)K_{\ell'}(t)h(u)dtdu - \int K_{\ell'}(\gamma - v)h(v)dv \notag \\
&= \int\int K_\ell(v-u)K_{\ell'}(\gamma - v)h(u)dudv - \int K_{\ell'}(\gamma - v)h(v)dv \notag \\
&= \int K_{\ell'}(\gamma - v)\left(\int K_\ell(v-u)h(u)du - h(v) \right)dv \notag \\
&= \int K_{\ell'}(\gamma - v)\left(K_\ell\star h(v) - h(v) \right)dv.  \notag
\end{align}
Hence, we derive
\begin{equation}\label{proof-th1-El}
\ch{\kv[\hllp] - \kv[\hlp]} = \ch{K_{\ell'}\star (K_\ell\star h - h)} \le \chuan{K}{1}\ch{K_\ell\star h - h},
\end{equation}
where the right hand side does not depend on $\ell'$ allowing us to take $\usup{\ell'\in H}$ in the left hand side. 

\par Thus, \eqref{eq:A1A2}, \eqref{eq:proof-th1-xi} and \eqref{proof-th1-El} give
\[
A_1 + A_2 \le 2\xi_T(\lb) + 2\chuan{K}{1}\ch{K_\ell\star h - h} + 2\frac{\chi\ch{K}}{\sqrt{M_T \ell}}.
\]
Then,
\begin{equation}\label{eq:proof-th1-A1A2}
\kv\left[(A_1 + A_2)^2\right] 
\le 12\kv[\xi_T^2(\lb)] + 12\chuan{K}{1}^2\ch{K_\ell\star h - h}^2 + 12\frac{\chi^2\ch{K}^2}{\ell}
\kv\left[\frac{1}{M_T}\right]. 
\end{equation}

\noindent For the term $A_3$, we have from \eqref{eq:upper-bound-var}
\begin{align}
\kv\left[A_3^2\right] &= \ch{\kv[\hl] - h}^2 + \kv\left[\ch{\hl - \kv[\hl]}^2\right] \notag\\
&\le \ch{K_\ell\star h - h}^2 + \frac{\ch{K}^2}{\ell}\kv\left[\frac{1}{M_T}\right]. \notag
\end{align}

\noindent Finally, replacing $\chi$ by $(1+\epsilon)(1+\chuan{K}{1})$, we have for any $\lb\in H$
\begin{align}
&\kv\left[\ch{\hat{h} - h}^2 \right] \le 2\kv\left[(A_1 + A_2)^2\right] + 2\kv\left[A_3^2\right] \notag\\
&\phantom{aaa} \le 24\kv\left[\xi^2_T(\lb)\right] + 2\left(1 + 12\chuan{K}{1}^2\right)\ch{K_\ell\star h - h}^2 \notag\\
&\phantom{aaaaaaaa} + 2\Big(1 +12(1 + \epsilon)^2(1 + \chuan{K}{1})^2\Big)\frac{\ch{K}^2}{\ell}\kv\left[\frac{1}{M_T}\right] \notag\\
&\phantom{aaa} \le 24\kv\left[\xi^2_T(\ell)\right] + C_1\left( \ch{K_\ell\star h - h}^2 + \frac{\ch{K}^2}{\ell}\kv\left[\frac{1}{M_T}\right]\right), \label{eq:proof-th1-bound1}
\end{align}
with $C_1$ a constant depending on $\epsilon$ and $\chuan{K}{1}$.

\par It remains to deal with the term $\kv\left[\xi_T^2(\ell)\right]$ where $\xi_T(\ell)$ is defined in \eqref{eq:proof-th1-xi}, 
\begin{align*}
\xi_T(\ell) &\le \usup{\ell'\in H}\left\{ \ch{\hllp - \kv[\hllp]} + \ch{\hlp - \kv[\hlp]}
	- \frac{\chi\ch{K}}{\sqrt{M_T\ell'}} \right\}_{+} \\
&\le \usup{\ell'\in H}\left\{\ch{\hlp - \kv[\hlp]}\chuan{K}{1}  + \ch{\hlp - \kv[\hlp]}
	- \frac{\chi\ch{K}}{\sqrt{M_T\ell'}} \right\}_{+} \\
&\le \usup{\ell'\in H}\left\{\big(1 + \chuan{K}{1}\big)\ch{\hlp - \kv[\hlp]}
	- \frac{(1+\epsilon)(1+\chuan{K}{1})\ch{K}}{\sqrt{M_T\ell'}} \right\}_{+} \\
&\le (1+\chuan{K}{1})S_T,
\end{align*}
where
\[
S_T := \usup{\ell\in H}\left\{\ch{\hl - \kv[\hl]} - \frac{(1+\epsilon)\ch{K}}{\sqrt{M_T\ell}}\right\}_+.
\]
Hence,
\begin{align}
\kv\big[\xi_T^2(\ell)\big] \le (1 + \chuan{K}{1})^2\kv\Big[\kv\big[S_T^2|M_T\big]\Big]. \label{eq:proof-th1-bound-xi}
\end{align}
If we show that
\begin{equation}\label{eq:proof-th1-bound0}
\kv\big[S_T^2|M_T = n\big] \le C_*\frac{1}{n},
\end{equation}
then
\begin{equation}\label{eq:proof-th1-bound2}
\kv\big[\xi_T^2(\ell)\big] \le C_*(1 + \chuan{K}{1})^2\kv\left[\frac{1}{M_T}\right]
\end{equation}
where $C_*$ is a constant.
\par Let us establish \eqref{eq:proof-th1-bound0}. When $M_T = n$, $\forall n\in\nhoa$, we set
\[
\kv\big[\Sigma^2_n\big] = \kv\big[S_T^2|M_T = n\big]
\] 
where 
\[
\Sigma_n := \usup{\ell\in H}\left\{\ch{Z_\ell} - \frac{(1+\epsilon)\ch{K}}{\sqrt{n\ell}}\right\}_+,
\]
with
\[
Z_\ell = \hl - \kv[\hl] = \frac{1}{n}\sum_{i=1}^n \Kl - \kv[\Kl].
\]
Then,
\begin{align*}
\kv\big[\Sigma^2_n\big] &= \kv\left[\usup{\ell\in H}\left\{\ch{Z_\ell} -  \frac{(1+\epsilon)\ch{K}}{\sqrt{n\ell}}\right\}_+^2 \right] \\
&\le \int_0^{+\vc} \xs\left[\usup{\ell\in H}\left\{\ch{Z_\ell} -  \frac{(1+\epsilon)\ch{K}}{\sqrt{n\ell}}\right\}_+^2 \ge x \right]dx \\
&\le \sum_{\ell\in H}\int_0^{+\vc} \xs\left[\left\{\ch{Z_\ell} -  \frac{(1+\epsilon)\ch{K}}{\sqrt{n\ell}}\right\}_+^2 \ge x \right]dx.
\end{align*}	
We bound this with Talagrand's inequality.

\par Let $\ahoa$ be a countable dense subset of the unit ball of $\mathbb{L}_2([0,1])$. We express the  norm $\ch{Z_\ell}$ as
\begin{align}
\ch{Z_\ell} &= \usup{a\in \ahoa}\int a(\gamma)Z_\ell(\gamma)d\gamma \notag \\
	&= \usup{a\in \ahoa}\sum_{i=1}^n \int a(\gamma)\frac 1n\left(\Kl - \kv[\Kl] \right)d\gamma.
	\notag
\end{align}
Let
\[
V_{i,\Gamma} = \int a(\gamma)\frac 1n\left(\Kl - \kv[\Kl] \right)d\gamma.
\]
Then $V_{i,\Gamma}$, $i=1,\ldots,n$ is a sequence of i.i.d random variables with zero mean. Thus, we can apply Talagrand's inequality 
(see \cite[p. 170]{Massart03}) to $\ch{Z_\ell} = \usup{a\in\ahoa}\sum_{i=1}^n V_{i,\Gamma}$. For 
all $\eta, x >0$, one has
\[
\xs\left(\ch{Z_\ell} \ge (1+\eta)\kv[\ch{Z_\ell}] + \sqrt{2\nu x} + c(\eta)bx \right) \le e^{-x},
\]
where $c(\eta) = 1/3 + \eta^{-1}$,
\[
\nu = \frac 1n \usup{a\in \ahoa} \kv\left[\left(\int a(\gamma)\left(\Kll -\kv[\Kll] \right)d\gamma \right)^2 \right], 
\]
and,
\[
b = \frac 1n \usup{y\in (0,1), a\in\ahoa}\int a(\gamma)\left(K_\ell(\gamma - y) - \kv[\Kll] \right)d\gamma.
\]

Next, we calculate the terms $\kv[\ch{Z_\ell}]$, $\nu$ and $b$. Applying Cauchy - Schwarz's inequality and using independence of variables, we get
\begin{align*}
\kv\big[\ch{Z_\ell}\big] &\le \left(\kv\big[\ch{Z_\ell}^2\big]\right)^{1/2} \\
&\le \left(\kv\left[\int \left(\frac 1n \sum_{i=1}^n \Kl - \kv[\Kl] \right)^2d\gamma \right] \right)^{1/2} 
\end{align*}
\begin{align*}
\phantom{\kv\big[\ch{Z_\ell}\big]} &= \frac 1n\left(\int \kv\left[\left(\sum_{i=1}^n \Kl - \kv[\Kl] \right)^2\right]d\gamma \right)^{1/2} \\
&= \frac 1n\left(\int \sum_{i=1}^n \kv\left[\left(\Kl - \kv[\Kl] \right)^2\right]d\gamma \right)^{1/2} \\
&\le \frac 1n\left(n\int \kv\left[K_\ell(\gamma - \Gamma_1^1)^2 \right]d\gamma \right)^{1/2} = \frac{\ch{K}}{\sqrt{n\ell}}.
\end{align*}
For the term $\nu$, we have
\begin{align*}
\nu &\le \frac 1n \usup{a\in \ahoa}\kv\left[\left(\int a(\gamma)K_\lb(\gamma - \Gamma_1^1)d\gamma \right)^2 \right] \\
&\le \frac 1n \usup{a\in \ahoa}\kv\left[\int |K_\ell(\gamma - \Gamma_1^1) |d\gamma
\times \int a^2(\gamma)|K_\ell(\gamma - \Gamma_1^1)|d\gamma \right] \\
&\le \frac{\chuan{K}{1}}{n}\usup{a\in \ahoa}\kv\left[\int a^2(\gamma)|K_\ell(\gamma - \Gamma_1^1)|d\gamma \right] \\
&\le \frac{\chuan{K}{1}}{n}\usup{a\in \ahoa}\int a^2(\gamma)\kv\left[|K_\ell(\gamma - \Gamma_1^1)| \right]d\gamma \\
&\le \frac{\chuan{K}{1}}{n}\usup{a\in \ahoa}\int\int a^2(\gamma)|K_\ell(\gamma - u)|h(u)dud\gamma \\
&\le \frac{\chuan{h}{\vc}\chuan{K}{1}^2}{n}.
\end{align*}
For the term $b$, we have
\begin{align*}
b &= \frac 1n \usup{y\in (0,1)}\ch{K_\lb(\cdot - y) - \kv[K_\lb(\cdot - \Gamma_1^1)]} \\
&\le \frac 1n \left(\usup{y\in (0,1)}\ch{K_\lb(\cdot - y)} + \left(\kv\big[\int K_\lb^2(\gamma - \Gamma_1^1)d\gamma\big] \right)^{1/2} \right) \le \frac{2\ch{K}}{n\sqrt{\lb}}.
\end{align*}
So, for all $\eta, x>0$, we have
\[
\xs\left(\ch{Z_\ell} \ge (1+\eta)\frac{\ch{K}}{\sqrt{n\ell}} + 
\chuan{h}{\vc}^{1/2}\chuan{K}{1}\sqrt{\frac{2x}{n}} + 2c(\eta)\frac{\ch{K}x}{n\sqrt{\ell}}
 \right) \le e^{-x}.
\]

Let $W_\lb$ be some strictly positive weights, we apply the previous inequality to $x=W_\lb + u$ for $u>0$. We have
\begin{align*}
\xs\Bigg(\ch{Z_\ell} \ge (1+\eta)\frac{\ch{K}}{\sqrt{n\ell}} &+ \chuan{h}{\vc}^{1/2}\chuan{K}{1}\sqrt{\frac{W_\lb}{n}} + 2c(\eta)\frac{\ch{K}W_\lb}{n\sqrt{\ell}} \\
&\phantom{aaaa} + \chuan{h}{\vc}^{1/2}\chuan{K}{1}\sqrt{\frac{u}{n}} + 2c(\eta)\frac{\ch{K}u}{n\sqrt{\ell}}\Bigg) \le e^{-W_\ell - u}.
\end{align*}
If we set
\[
\Psi_\lb = (1+\eta)\frac{\ch{K}}{\sqrt{n\ell}} + \chuan{h}{\vc}^{1/2}\chuan{K}{1}\sqrt{\frac{W_\lb}{n}} + 2c(\eta)\frac{\ch{K}W_\lb}{n\sqrt{\ell}},
\]
then,
\[
\xs\Big(\ch{Z_\ell} - \Psi_\ell \ge \chuan{h}{\vc}^{1/2}\chuan{K}{1}\sqrt{\frac{u}{n}} + 2c(\eta)\frac{\ch{K}u}{n\sqrt{\ell}} \Big)\le e^{-W_\ell - u}.
\]
Let
\[
\Lambda = \kv\left[\usup{\lb\in H}\left(\ch{Z_\ell} - \Psi_\lb \right)^2_+ \right] =\int_0^{+\vc} \xs\left[\usup{\lb\in H}\left(\ch{Z_\ell} - \Psi_\lb \right)^2_+ \ge x \right]dx .
\]
An upper bound of $\Lambda$ is given by
\begin{align*}
\Lambda &\le \sum_{\lb\in H}\int_0^{+\vc} \xs\left[\left(\ch{Z} - \Psi_\lb \right)^2_+ \ge x \right]dx.
\end{align*}
Let us take $u$ such that 
\[
x = f(u)^2 = \left(\chuan{h}{\vc}^{1/2}\chuan{K}{1}\sqrt{\frac{u}{n}} + 2c(\eta)\frac{\ch{K}u}{n\sqrt{\ell}} \right)^2.
\]
So,
\[
dx = 2f(u)\left(\chuan{h}{\vc}^{1/2}\chuan{K}{1}\frac{1}{2\sqrt{nu}} +  2c(\eta)\frac{\ch{K}}{n\sqrt{\ell}} \right)du.
\]
Hence,
\begin{align}
\Lambda &\le \sum_{\lb\in H}\int_0^{+\vc} e^{-W_\lb - u}2f(u)\left(\chuan{h}{\vc}^{1/2}\chuan{K}{1}\frac{1}{2\sqrt{nu}} +  2c(\eta)\frac{\ch{K}}{n\sqrt{\ell}} \right)du \notag\\
&\le \sum_{\lb\in H}\int_0^{+\vc} e^{-W_\lb - u}2f(u)\left(\chuan{h}{\vc}^{1/2}\chuan{K}{1}\sqrt{\frac{u}{n}} + 2c(\eta)\frac{\ch{K}u}{n\sqrt{\ell}} \right)u^{-1}du \notag\\
&\le 2\sum_{\lb\in H}e^{-W_\lb}\int_0^{+\vc}f^2(u)e^{-u}u^{-1}du  \notag\\
&\le C_{\eta}\sum_{\lb\in H}e^{-W_\lb}\left(\chuan{h}{\vc}\chuan{K}{1}^2\int_0^{+\vc}e^{-u}du 
+ \frac{\ch{K}^2}{\ell^2}\int_0^{+\vc}ue^{-u}du \right)\times \frac 1n \notag\\
&\le C_{\eta}\sum_{\lb\in H}e^{-W_\lb}\left(\chuan{h}{\vc}\chuan{K}{1}^2 
+ \frac{\ch{K}^2}{\ell^2} \right)\times \frac 1n. \label{eq:proof-th1-bound-R}
\end{align}

\noindent We need to choose $W_\ell$ and $\eta$ such that 
\begin{equation}
\kv\big[\Sigma^2_n \big] = \kv\left[\usup{\ell\in H}\left\{\ch{Z_\ell} - \frac{(1+\epsilon)\ch{K}}{\sqrt{n\ell}}\right\}_+^2 \right] \le \Lambda.
\end{equation}

\noindent Let $\theta >0$, we choose
\[
W_\lb = \frac{\theta^2\ch{K}^2}{2\chuan{h}{\vc}\chuan{K}{1}^2\sqrt{\ell}},
\]
the we have
\[
\Psi_\lb = (1+\eta)\frac{\ch{K}}{\sqrt{n\lb}} + \frac{\theta\ch{K}}{\sqrt{2n\sqrt{\ell}}} + 
\frac{c(\eta)\theta^2\ch{K}^3}{\chuan{h}{\vc}\chuan{K}{1}^2}\frac{1}{n\ell}.
\]
Obviously, the series in \eqref{eq:proof-th1-bound-R} is finite and for any $\ell \in H$, since $\lb \le 1$,
we have
\begin{align*}
\Psi_\lb &\le (1 + \eta + \theta)\frac{\ch{K}}{\sqrt{n\lb}} + 
\frac{c(\eta)\theta^2\ch{K}^3}{\chuan{h}{\vc}\chuan{K}{1}^2}\frac{1}{n\ell} \\
&\le \left(1 +\eta + \theta +  \frac{c(\eta)\theta^2\ch{K}^2}{\chuan{h}{\vc}\chuan{K}{1}^2}\frac{1}{\sqrt{n\ell}}\right)\frac{\ch{K}}{\sqrt{n\lb}}.
\end{align*}
Since $ H\subset \left\{\triangle^{-1}, \triangle=1,\ldots, \triangle_{\max} \right\}$, if we choose $\triangle_{\max} = \lfloor \delta n \rfloor $ for some $\delta >0$, then $\lb_{\min} = \triangle_{\max}^{-1}$ and we obtain
\[
\Psi_\lb \le  \left(1 +\eta + \theta +  \frac{c(\eta)\theta^2\ch{K}^2\sqrt{\delta}}{\chuan{h}{\vc}\chuan{K}{1}^2}\right)\frac{\ch{K}}{\sqrt{n\lb}}.
\]
It remains to choose $\eta = \epsilon/2$ and $\theta$ small enough such that
\[
\theta +  \frac{c(\eta)\theta^2\ch{K}^2\sqrt{\delta}}{\chuan{h}{\vc}\chuan{K}{1}^2} = \frac{\epsilon}{2},
\]
then 
\[
\Psi_\lb \le (1+\epsilon)\frac{\ch{K}}{\sqrt{n\lb}},
\]
and we get 
\[
\kv\left[\Sigma^2_n \right] \le C_*\times\frac{1}{n},
\]
where $C_*$ is a constant depending on $\delta$,$\epsilon$,$\chuan{h}{\vc}$,$\chuan{K}{1}$ and $\ch{K}$. Hence, we get \eqref{eq:proof-th1-bound2}.

\par Combining \eqref{eq:proof-th1-bound1} and \eqref{eq:proof-th1-bound2}, we obtain
\[
\kv\left[\ch{\hat{h} - h}^2 \right] \le C_1\left( \ch{K_\ell\star h - h}^2 + \frac{\ch{K}^2}{\ell}\kv\left[\frac{1}{M_T}\right]\right) + C_*\kv\left[\frac{1}{M_T}\right]
\]

Moreover, since $N_T > N_0$, we have
\begin{align}
\kv\left[\frac{1}{M_T}\right] &= \kv\left[\frac{1}{N_T-N_0}\right] =\kv\left[\frac{N_T}{N_T - N_0}\frac{1}{N_T} \right] =\kv\left[\frac{1}{1-\frac{N_0}{N_T}}\frac{1}{N_T} \right] \notag\\
&\le \kv\left[\frac{1}{1-\frac{N_0}{N_0+1}}\frac{1}{N_T} \right] \notag \\ 
&\le (N_0+1)\kv\left[\frac{1}{N_T}\right]. \label{eq:inverse-Mt}
\end{align}

Then, using \eqref{eq:inverse-N_T-geom}, \eqref{eq:bound-inverse-N_T-nb} and \eqref{eq:inverse-Mt}, recall the definition of $\varrho(T)^{-1}$ in \eqref{eq:rate}, we obtain for any $\lb \in  H$
\[
\kv\left[\ch{\hat h - h}^2 \right]\le C_1\left(\ch{K_\lb\star h - h}^2  + \frac{\ch{K}^2}{\lb}\varrho(T)^{-1} \right) + C_2\varrho(T)^{-1} .
\]
This ends the proof of Theorem \ref{th1}.

\subsection{Proof of Theorem \ref{th:upper-bound}}
We begin with the bias term $\ch{K_\ell\star h - h}$ in the right hand side of the oracle inequality \eqref{eq:oracle-ineq}. For any $\lb\in H$ and $\gamma\in (0,1)$, let $k = \lfloor \beta \rfloor$ and $b(\gamma) = K_\ell\star h(\gamma) - h(\gamma)$, then we have
\[
h(\gamma + u\ell) = h(\gamma) + h'(\gamma)u\ell + \cdots + \frac{(u\ell)^k}{(k-1)!}\int_0^1 (1-\theta)^{k-1}h^{(k)}(\gamma + 
\theta u\ell)d\theta. 
\]
Since $K$ is a kernel of order $\beta^*$ and $\beta\in (0,\beta^*)$, we get
\begin{align*}
b(\gamma) =\int K(u)\frac{(u\ell)^k}{(k-1)!}\left[\int_0^1 (1-\theta)^{k-1}\left(h^{(k)}(\gamma + \theta u\ell) - h^{(k)}(\gamma)\right)d\theta\right]du.
\end{align*}
Setting $E_{k,\ell}(u) = |K(u)|\frac{|u\ell|^k}{(k-1)!}$ for the sake of notation. Since $h\in \mathcal{H}(\beta, L)$ and applying twice the generalized Minskowki's inequality, we obtain
\begin{align*}
&\ch{h - \kv[\hat h]}^2 = \int b^2(\gamma)d\gamma \\
&\le \int\left(\int E_{k,\ell}(u)\left[\int_0^1 (1-\theta)^{k-1}\big|h^{(k)}(\gamma + \theta u\ell) - 
      h^{(k)}(\gamma)\big|d\theta\right]du \right)^2d\gamma \\
&\le \Bigg(\int E_{k,\ell}(u)\left[\int\left(\int_0^1 (1-\theta)^{k-1}\big|h^{(k)}(\gamma + \theta u\ell) - 
      h^{(k)}(\gamma)\big|d\theta \right)^2d\gamma \right]^{1/2}du \Bigg)^2 \\
&\le \Bigg(\int E_{k,\ell}(u) \left[\int_0^1 (1-\theta)^{k-1}\left(\int \big|h^{(k)}(\gamma + \theta u\ell) - 
      h^{(k)}(\gamma)\big|^2 d\gamma \right)^{1/2}d\theta \right]du \Bigg)^2 \\
&\le \left(\int E_{k,\ell}(u) \left[\int_0^1
(1-\theta)^{k-1} L(\theta u\ell)^{\beta-k} d\theta \right]du \right)^2 \\
&\le \left(\int |K(u)|\frac{|u\ell|^k}{(k-1)!} \left[\int_0^1
(1-\theta)^{k-1}L(u\ell)^{\beta - k}d\theta \right]du \right)^2 \\
&\le  C_{K,L,\beta}\ell^{2\beta},
\end{align*}
where $C_{K,L,\beta} = \left(\frac{L}{k!}\int |u|^\beta|K(u)|du \right)^2$.

\par Finally, we have
\begin{equation}\label{eq:temp3}
\kv\left[\ch{\hat h - h}^2 \right]\le C_1\underset{\ell\in H}{\inf}\left\{C_{K,L,\beta}\ell^{2\beta}  + \frac{\ch{K}^2}{\lb}\varrho(T)^{-1} \right\} + C_2\varrho(T)^{-1} .
\end{equation}
Taking the derivative of the expression inside the $\underset{\lb\in H}{\inf}$ of \eqref{eq:temp3} with respect to $\lb$, we obtain the minimizer
\[
\ell^* = \left(\frac{\ch{K}^2}{2\beta C_{K,L,\beta}} \right)^{\frac{1}{2\beta+1}}\varrho(T)^{-\frac{1}{2\beta+1}}. 
\]
Since the optimal bandwidth $\hat\lb$ is proportional to $\ell^*$ up to a multiplicative constant. Therefore, by substituting $\ell$ by $\hat\lb$ in the right hand side of \eqref{eq:temp3}, we obtain
\[
\kv\left[\ch{\hat h - h}^2 \right]\le C_3\varrho(T)^{-\frac{2\beta}{2\beta+1}},
\]
with $C_3$ a constant depending on $N_0$, $\delta$, $\epsilon$, $\chuan{K}{1}$, $\ch{K}$, $\chuan{h}{\vc}$, $\beta$ and $L$.
This ends the proof of Theorem \ref{th:upper-bound}.

\subsection{Proof of Theorem \ref{th:lower-bound}}
For $T>0$, let us denote by $\hat{h}_T$ the estimator of $h$. To prove the Theorem \ref{th:lower-bound}, we apply the general reduction scheme proposed by Tsybakov \cite{Tsybakov04} (Section 2.2, p.79). We will show the existence of a family $\hhoa_{m,T} = \big\{h_{j,T}: j = 0,1,\ldots, m \big\}$ such that:

\begin{itemize}
\item[1)] $h_{j,T} \in \hhoa(\beta, L)$, $j=0,\ldots,m$. 

\item[2)] $\ch{h_{j,T} - h_{k,T}} \ge 2c\,e^{-\frac{\beta}{2\beta+1}RT},\; 0\le j < k \le m$.

\item[3)] $\dfrac{1}{m}\sum_{j=1}^m K(P_j, P_0) \le \vartheta\log(m)$ for $0 < \vartheta < 1/8$. $P_j$ and $P_0$ are the distribution of observations when the division kernels are $h_{j,T}$ and $h_0$, respectively. $K(P,Q)$ denotes the Kullback-Leibler divergence between two measures $P$ and $Q$:
\[
K(P,Q) = \begin{cases}
\int \log\frac{dP}{dQ}dP, &\text{ if } P\ll Q \\
+\vc, & \text{ otherwise}.
\end{cases}
\]
\end{itemize}
Under the preceding conditions 1, 2, 3, Tsybakov \cite{Tsybakov04} (Theorem 2.5, p.99) show that
\begin{equation}
\underset{\hat{h}_T}{\inf}\,\underset{h\in \hhoa_{m,T}}{\max}\,\xs\left(\ch{\hat{h}_T - h}^2 \ge c^2 e^{-\frac{2\beta}{2\beta + 1}RT} \right) \ge C',
\end{equation}
where the infimum is taken over all estimators $\hat{h}_T$ and positive constant $C'$ is independent of $T$. This will be sufficient to obtain Theorem \ref{th:lower-bound} by \cite[Theorem 2.7]{Tsybakov04}. The proof ends with proposing a family $\hhoa_{m,T}$ and checking the assumptions 1, 2, 3. \\[6pt]
\noindent\textit{\bf Construction of the family} $\hhoa_{m,T}$: 

\par The idea is construct a family of perturbations around $h_0$ which is a symmetric density with respect to $\frac 12$ and belongs to $\hhoa(\frac L2, \beta)$. For the simplification, we choose $h_0(\gamma) = \id_{(0,1)}(\gamma)$.

\par Let $c_0>0$ be a real number, and let $\gamma \in (0,1)$, $f(\gamma) = LD^{-\beta} g\left(D\gamma\right)$ where $g$ is a regular function having support $(0,1)$ and $\int g(\gamma)d\gamma = 0$, $g\in \hhoa(\frac 12,\beta)$, we define
\[
D = \lceil c_0 e^{\frac{RT}{2\beta + 1}}\rceil\;\text{ and }\; f_k(\gamma) = f\left(\gamma - \frac{(k-1)}{D} \right),
\]
By definition, the functions $f_k$'s have disjoint support and one can check that the functions $f_k \in \hhoa(\frac{L}{2}, \beta)$.

\par Then, the function $h_{j,T}$ will be chosen in 
\[
\mathcal{D} = \left\{h_\delta(\gamma) = h_0(\gamma) + c_1\sum_{k=1}^D \delta_k f_k(\gamma): \delta = (\delta_1,\ldots, \delta_D) \in \{0,1\}^D \right\},
\]
where 
\begin{equation}\label{eq:cond-positive-density}
c_1 = \min\left(\frac{1}{LD^{-\beta}\chuan{g}{\vc}}, 1 \right).
\end{equation}

We now check that $h_\delta$ is a density, since $\int h_\delta(\gamma)d\gamma = \int h_0(\gamma)d\gamma = 1$, it remains to verify that $h_\delta(\gamma) \ge 0 \;\forall\, \gamma$. We have 
\begin{align*}
\underset{(0,1)}{\inf} h_\delta(\gamma) &\ge \underset{(0,1)}{\inf}h_0 - \chuan{c_1\sum_{k=1}^D \delta_kf_k}{\vc} \\
&\ge 1 - c_1 LD^{-\beta}\underset{k}{\max}\,\underset{\gamma}{\sup}\,|\delta_k|g\big(D\gamma - (k-1) \big)\\
&\ge 1 - c_1 LD^{-\beta}\chuan{g}{\vc} \ge 0,
\end{align*}
by the choice of $c_1$. Thus the family of densities $\mathcal{D}$ is well-defined.\\[6pt]
\noindent\textit{\bf 1) The condition $h_{j,T}\in \hhoa(\beta, L)$}:\\[6pt]
Let us denote $q = \lfloor \beta \rfloor$, then for all $\gamma, \gamma' \in (0,1)$ we have
\begin{align*}
\Big| h_\delta^{(q)}(\gamma) - h_\delta^{(q)}(\gamma')\Big| &=
\Big|h_0^{(q)}(\gamma) - h_0^{(q)}(\gamma') + c_1\sum_{k=1}^D \delta_k\big(f_k^{(q)}(\gamma) - f_k^{(q)}(\gamma') \big) \Big| \\
&\le c_1\sum_{k=1}^D|\delta_k|\Big|f_k^{(q)}(\gamma) - f_k^{(q)}(\gamma') \Big| \\
&\le c_1\underset{k}{\max}\Big|f_k^{(q)}(\gamma) - f_k^{(q)}(\gamma') \Big| \\
&\le c_1 LD^{-\beta}\underset{k}{\max}\,D^{q}\Big|g^{(q)}(D\gamma - (k-1)) - g^{(q)}(D\gamma' - (k-1))\Big| \\
&\le c_1 LD^{\lfloor \beta \rfloor - \beta}D^{\beta - \lfloor \beta \rfloor}|\gamma - \gamma'|^{\beta - \lfloor \beta \rfloor} \le L|\gamma - \gamma'|^{\beta - \lfloor \beta \rfloor},
\end{align*}
which is always satisfied with $c_1 = \min\left(\frac{1}{LD^{-\beta}\chuan{g}{\vc}}, 1 \right)$, thus $h_\delta \in \hhoa(L,\beta)$.\\[6pt]
\noindent\textit{\bf 2) The condition $\ch{h_{j,T} - h_{k,T}} \ge 2c\,e^{-\frac{\beta}{2\beta+1}RT}$}: \\[4pt]
For all $\delta, \delta' \in \{0,1\}^D$, we have
\begin{align*}
\ch{h_\delta - h_{\delta'}} &= \left[\int_0^1\big(h_\delta(\gamma) - h_{\delta'}(\gamma)\big)^2d\gamma \right]^{1/2} = \left[\int_0^1\left(c_1\sum_{k=1}^D (\delta_k - \delta_k')f_k(\gamma) \right)^2d\gamma \right]^{1/2} \\
&= c_1\left[\int_0^1\sum_{k=1}^D (\delta_k - \delta_k')^2f^2_k(\gamma) d\gamma \right]^{1/2} = c_1\left[\sum_{k=1}^D (\delta_k - \delta_k')^2\int_{\frac{k-1}{D}}^{\frac{k}{D}} f_k^2(\gamma)d\gamma\right]^{1/2}\\
&= c_1\left[\sum_{k=1}^D (\delta_k - \delta_k')^2\int_{\frac{k-1}{D}}^{\frac{k}{D}} L^2 D^{-2\beta} g^2\left(D\gamma - (k-1)\right)d\gamma\right]^{1/2} \\
&= c_1L D^{-\beta - 1/2}\ch{g}\left[\sum_{k=1}^D (\delta_k - \delta_k')^2 \right]^{1/2} =
c_1 L D^{-\beta - 1/2}\ch{g}\sqrt{d_H(\delta, \delta')},
\end{align*}
where $d_H(\delta, \delta') = \sum_{k=1}^D \id\{\delta_k \neq \delta_k'\}$ is the Hamming distance between $\delta$ and $\delta'$.

\par According to the Lemma of Varshamov-Gilbert (cf. Tsybakov \cite{Tsybakov04}, p.104), there exist a subset $\left\{\delta^{(0)},\ldots, \delta^{(m)} \right\}$ of $\{0,1\}^D$ with cardinal \eqref{eq:Varshamov-bound2} such that $\delta^{(0)} = (0,\ldots,0)$, 
\begin{equation}\label{eq:Varshamov-bound2}
m \ge 2^{D/8},
\end{equation} 
and
\begin{equation}\label{eq:Varshamov-bound1}
d_H(\delta^{(j)}, \delta^{(k)}) \ge \frac{D}{8}, \;\forall\; 0 \le j < k \le m.
\end{equation}
Then, by setting $h_{j,T}(x) = h_{\delta^{(j)}}(x)$, $j=0,\ldots, m$, we obtain
\begin{align*}
\ch{h_{j,T} - h_{k,T}} &= c_1 LD^{-\beta -\frac 12}\ch{g}\sqrt{d_H(\delta^{(j)}, \delta^{(k)})} \\
&\ge c_1 LD^{-\beta - 1/2}\ch{g}\sqrt{\frac{D}{8}}\\
&\ge \frac{c_1 L}{4}\ch{g}D^{-\beta}
\end{align*}
whenever $D\ge 8$.

\par Suppose that $N_T\ge N_{T^*}$ where $T^* = \log\left(\frac{7}{c_0} \right)\frac{2\beta+1}{R}$. Then, $D\ge 8$ and $D^\beta \le (2c_0)^\beta e^{\frac{\beta}{2\beta + 1} RT}$.
This implies: 
\[
\ch{h_{j,T} - h_{k,T}} \ge \frac{c_1 L}{4}\ch{g}(2c_0)^{-\beta} e^{-\frac{\beta}{2\beta+1}RT},
\]
But,
\[
\min\left(\frac{1}{L\chuan{g}{\vc}},1 \right) \le c_1 \le 1
\]
Hence, we obtain
\[
\ch{h_{j,T} - h_{k,T}} \ge 2c\;e^{-\frac{\beta}{2\beta+1}RT},
\]
where
\[
c = \frac{\min(1,L\ch{g})}{8}(2c_0)^{-\beta}.
\]
\noindent\textit{\bf 3) The condition $\frac{1}{m}\sum_{j=1}^m K(P_j, P_0) \le \vartheta\log(m)$ for $0<\vartheta<1/8$}:\\[6pt]
We need to show that for all $\delta \in \{0,1\}^D$,
\[
K(P_\delta, P_0) \le \vartheta \log(m),
\]
where 
\[
K(P_\delta, P_0) = \kv\left[\log\frac{dP_\delta}{dP_0}|_{\Ftr_T}(Z) \right],
\]
and where $(Z_t)_{t\in [0,T]}$ is defined in \eqref{eq:SDE-model} with the random measure $Q$ having intensity $q(ds,di,d\gamma) = Rh_\delta(\gamma)ds\,n(di)d\gamma$.

\par Here, the difficulty comes from the fact that $N_T$ is variable because the observations result from a stochastic process $Z_t$. The law of these observations is not a probability distribution on a fixed $\rhoa^n$ where $n$ would be the sample size, but rather a probability distribution on a path space. $P_\delta$ is the probability distribution when the Poisson point measure $Q$ has intensity $Rh_\delta(\gamma)ds\,n(di)d\gamma$. Thus a natural tool is to use Girsanov's theorem (see \cite{Jacob87}, Theorem 3.24, p. 159) saying that $P_\delta$ is absolutely continuous with respect to $P_0$ on $\Ftr_T$ with
\[
\frac{dP_\delta}{dP_0}|_{\Ftr_T} = \mathfrak{D}^\delta_T,
\]
where $(\mathfrak{D}^\delta_t)_{t\in [0,T]}$ is the unique solution of the following SDE (see Proposition 4.17 of \cite{Tran06} for a similar SDE):
\begin{equation}\label{eq:SDE-exp-martingale}
 \mathfrak{D}^\delta_T = 1 + \int_0^T\int_{\ehoa}  \mathfrak{D}^\delta_{s-}\id_{\{i\le N_{s-}\}}\left(\frac{h_\delta(\gamma)}{h_0(\gamma)} - 1\right)Q(ds,di,d\gamma).
\end{equation}
Apply It\^{o} formula for jump processes to \eqref{eq:SDE-exp-martingale}, we get 
\begin{align*}
\log \mathfrak{D}^\delta_T &= \int_0^T\int_{\ehoa} \id_{\{i\le N_{s-}\}}\Bigg[\log\Bigg(\mathfrak{D}^\delta_{s-} - \bigg(\frac{h_\delta(\gamma)}{h_0(\gamma)} - 1\bigg)\mathfrak{D}^\delta_{s-} \Bigg) - \log \mathfrak{D}^\delta_{s-} \Bigg]Q(ds,di,d\gamma) \\
&= \int_0^T\int_{\ehoa} \id_{\{i\le N_{s-}\}}\log\frac{h_\delta(\gamma)}{h_0(\gamma)} Q(ds,di,d\gamma)
= \sum_{i=1}^{N_T}\log\frac{h_\delta(\Gamma^1_i)}{h_0(\Gamma^1_i)}
\end{align*}
by definition of $(\Gamma^1_1, \ldots, \Gamma^1_{N_T})$.\\[3pt]
Then,
\begin{align*}
K(P_\delta, P_0) &= \kv_\delta\left[\log \mathfrak{D}^\delta_T \right] = \kv_\delta\left[\sum_{i=1}^{N_T}\log\frac{h_\delta(\Gamma^1_i)}{h_0(\Gamma^1_i)} \right]\\
&= \kv\left[N_T \right]\kv_\delta\left[\log\frac{h_\delta(\Gamma^1_1)}{h_0(\Gamma^1_1)} \right]
= \kv\left[N_T \right]\int_0^1 h_\delta(\gamma)\log\frac{h_\delta(\gamma)}{h_0(\gamma)}d\gamma.
\end{align*}
Here, $\kv\left[N_T\right]$ does not depend on $h_\delta$ and we have $\kv[N_T] = N_0e^{RT}$.
Thus, recall the definition of $h_\delta(\cdot)$ and note that $\log(1+x)\le x$ for $x > -1$, we get 
\begin{align*}
K(P_\delta, P_0) &= N_0e^{RT}\int_0^1 h_{\delta}(\gamma)\log(h_\delta (\gamma))d\gamma \\
&= N_0e^{RT} \int_0^1 \Big(1 + c_1\sum_{k=1}^D \delta_k f_k(\gamma)\Big)\log\Big(1 + c_1\sum_{k=1}^D \delta_k f_k(\gamma) \Big)d\gamma \\
&= N_0e^{RT} \sum_{k=1}^D \int_{\frac{k-1}{D}}^{\frac{k}{D}} \big(1 + c_1\delta_k f_k(\gamma)\big)\log\big(1 + c_1\delta_kf_k(\gamma)\big)d\gamma \\
&= N_0e^{RT}\sum_{k=1}^D \delta_k \int_0^{1/D} \big(1 + c_1f(\gamma)\big)\log\big(1 + c_1f(\gamma)\big)d\gamma \\
&\le N_0e^{RT}D \int_0^{1/D} \big(1 + c_1f(\gamma)\big)c_1 f(\gamma)d\gamma \\
&\le N_0e^{RT}\left[c_1LD^{-\beta}\int_0^{1/D}g(D\gamma)Dd\gamma +  c_1^2L^2D^{-2\beta}\int_0^{1/D}g^2(D\gamma)Dd\gamma\right] \\
&\le N_0e^{RT}c_1^2L^2D^{-2\beta}\int_0^1 g^2(\gamma)d\gamma  \\
&\le N_0 c_1^2L^2\ch{g}^2 e^{RT} c_0^{-2\beta} e^{-\frac{2\beta}{2\beta+1}RT} \\
&\le N_0 L^2\ch{g}^2 c_0^{-2\beta-1}D \quad \text{ since } c_1 \le 1.
\end{align*}
From \eqref{eq:Varshamov-bound2}, we have $m \ge 2^{D/8}$ then 
\[
D \le \frac{8\log(m)}{\log(2)}.
\]
Hence, if we set
\[
c_0 = \left(\frac{8 N_0 L^2\ch{g}^2}{\vartheta \log(2)} \right)^{1/(2\beta + 1)},
\]
we obtain $K(P_\delta, P_0) \le \vartheta \log(m)$. This ends the proof of Theorem \ref{th:lower-bound}.\newpage

\noindent \textbf{Acknowledgements}\\[6pt]
The author is deeply grateful to V.C. Tran, V. Rivoirard and T.M. Pham Ngoc for the guidance and useful suggestions. The author would like to express his sincere thanks to P. Massart for constructive comments that improved the results. The author would like to thank Centre de Resources Informatiques de Lille 1 (CRI) for the computational cluster to implement the numerical simulations. This work is supported by Program 911 of Vietnam Ministry of Education and Training and is partly supported by the French Agence Nationale de la Recherche (ANR 2011 BS01 010 01 projet Calibration).

\fontsize{10}{9.6}
\selectfont

\end{document}